\documentclass[11pt]{amsart}
\usepackage[top=1.5in, bottom=1.5in, left=1.5in, right=1.5in]{geometry}
\expandafter\let\csname ver@amsthm.sty\endcsname\relax

\usepackage[all]{xy}
\usepackage{graphicx} 
\usepackage{tikz}
\usetikzlibrary{patterns}
\usetikzlibrary{decorations.pathreplacing}
\usetikzlibrary{arrows}
\usetikzlibrary{decorations.markings}

\usepackage{amsmath}
\usepackage{amssymb}
\usepackage{enumerate}
\usepackage{mathdots}
\usepackage{mathtools}

\usepackage{hyperref}
\usepackage{amsthm}
\usepackage[capitalize,noabbrev]{cleveref}

\allowdisplaybreaks

\usepackage{ytableau}

\numberwithin{equation}{section}

\newtheorem{thm}{Theorem}[section]
\newtheorem{lemma}[thm]{Lemma}
\newtheorem{cor}[thm]{Corollary}
\newtheorem{prop}[thm]{Proposition}

\newtheorem{Definition}[thm]{Definition}
\newenvironment{definition}
  {\begin{Definition}\rm}{\end{Definition}}

\newtheorem{Example}[thm]{Example}
\newenvironment{example}
  {\begin{Example}\rm}{\hfill$\lozenge$\end{Example}}

\newtheorem{Remark}[thm]{Remark}
\newenvironment{remark}
  {\begin{Remark}\rm}{\hfill$\lozenge$\end{Remark}}
  
\crefname{thm}{Theorem}{Theorems}
\crefname{lemma}{Lemma}{Lemmas}
\crefname{cor}{Corollary}{Corollaries}
\crefname{prop}{Proposition}{Propositions}
\crefname{conj}{Conjecture}{Conjectures}
\crefname{question}{Question}{Questions}
\crefname{problem}{Problem}{Problems}

\crefname{definition}{Definition}{Definitions}
\crefname{example}{Example}{Examples}
\crefname{remark}{Remark}{Remarks}

\newcommand{\N}{\mathbb{N}}

\newcommand{\Q}{\mathbb{Q}}
\newcommand{\R}{\mathbb{R}}

\newcommand{\E}{\mathbb{E}}

\newcommand{\SYT}{\mathcal{SYT}}
\newcommand{\SYTo}{\SYT^{+1}}
\newcommand{\PP}{\mathcal{PP}}
\newcommand{\PPo}{\PP^{+1}}
\newcommand{\J}{\mathcal{J}}
\newcommand{\Lin}{\mathcal{L}}
\newcommand{\Lino}{\Lin^{+1}}
\newcommand{\RPP}{\mathcal{RPP}}
\newcommand{\RPPo}{\RPP^{+1}}

\DeclareMathOperator{\Des}{D}
\DeclareMathOperator{\comaj}{comaj}
\DeclareMathOperator{\maj}{maj}
\DeclareMathOperator{\rk}{rk}
\DeclareMathOperator{\ddeg}{ddeg}

\newcommand{\Tin}[1]{\mathcal{T}_{#1}^+}
\newcommand{\Tout}[1]{\mathcal{T}_{#1}^-}
\newcommand{\Tq}[2]{\mathcal{T}^{#2}_{#1}}

\newcommand{\muuni}[1]{\mu_{\mathrm{uni}}^{#1}}
\newcommand{\murpp}[2]{\mu_{\leq,#2}^{#1}}
\newcommand{\mulin}[1]{\mu_{\mathrm{lin}}^{#1}}
\newcommand{\murk}[1]{\mu_{\rk}^{#1}}

\newcommand{\shape}[1]{P_{#1}}
\newcommand{\sshape}[1]{P^{\mathrm{shift}}_{#1}}

\newcommand{\rect}[2]{#1 \times #2}
\newcommand{\sstair}[1]{\delta_{#1}}

\newcommand{\bracket}[1]{\chi(#1)}

\newcommand{\wt}{\vartheta}

\newcommand{\ustar}{u_*}
\newcommand{\pstar}{p_*}
\newcommand{\istar}{i_*}
\newcommand{\rstar}{r_*}
\newcommand{\dstar}{d_*}

\DeclareRobustCommand{\qbinom}{\genfrac{\lbrack}{\rbrack}{0pt}{}}

\newcommand{\emailhref}[1]{\email{\href{#1}{#1}}}
\newcommand{\dfn}[1]{\textcolor{blue}{\emph{#1}}}

\title[\MakeLowercase{\texorpdfstring{$q$}{q}}-enumeration of barely set-valued tableaux and plane partitions]{On the \MakeLowercase{\texorpdfstring{$q$}{q}}-enumeration of barely set-valued tableaux and plane partitions}

\author[S. Hopkins]{Sam Hopkins}\emailhref{samuelfhopkins@gmail.com}
\address{Department of Mathematics, Howard University, Washington, DC 20059, USA}

\author[A. Lazar]{Alexander Lazar}\emailhref{alexander.leo.lazar@ulb.be}
\address{Department of Mathematics, Universit\'{e} Libre de Bruxelles, Brussels, Belgium}

\author[S. Linusson]{Svante Linusson}\emailhref{linusson@math.kth.se}
\address{Department of Mathematics, KTH Royal Institute of Technology, Stockholm, Sweden}


\date{\today}

\begin{document}

\begin{abstract}
Barely set-valued tableaux are a variant of Young tableaux in which one box contains two numbers as its entry. It has recently been discovered that there are product formulas enumerating certain classes of barely set-valued tableaux. We give some $q$-analogs of these product formulas by introducing a version of major index for these tableaux. We also give product formulas and $q$-analogs for barely set-valued plane partitions. Many of the results are stated in the generality of $P$-partitions that then specialize to particularly nice formulas for rectangles and minuscule posets. The proofs use several probability distributions on the set of order ideals of a poset, depending on the real parameter $q>0$, which we think could be of independent interest.
\end{abstract}

\maketitle

\section{Introduction and statement of results} \label{sec:intro}

Let $\lambda$ be a partition of $n$. Recall that a \dfn{standard Young tableau (SYT)} of shape~$\lambda$ is a filling of the boxes of the Young diagram of $\lambda$ with the numbers $1,\ldots,n$, each appearing once, so that entries are increasing along rows and down columns. Standard Young tableaux are central objects in algebraic and enumerative combinatorics. We use $\SYT(\lambda)$ to denote the set of SYTs of shape~$\lambda$. The celebrated \dfn{hook length formula} of Frame, Robinson, and Thrall~\cite{frame1954hook} asserts that
\begin{equation} \label{eqn:hlf}
\#\SYT(\lambda) = n! \cdot \prod_{u \in \lambda} \frac{1}{h(u)}
\end{equation}
where $h(u)$ denotes the \dfn{hook length} of the box $u$ of $\lambda$. See, e.g., \cite{sagan1990ubiquitous} for a modern discussion of the hook length formula.

In this paper we study \emph{set-valued} tableaux, with an eye towards proving product formulas analogous to~\eqref{eqn:hlf} in some special cases. A \dfn{set-valued tableau} is like a usual tableau except that its entries are finite, non-empty sets of numbers, rather than single numbers. More precisely, a \dfn{standard set-valued tableau} of shape $\lambda$ is a set-valued filling of the  boxes of the Young diagram of $\lambda$ with the numbers $1,\ldots,n+k$ for some $k\geq 0$, with each number appearing once, so that entries are increasing along rows and down column in the sense that $\max(u) < \min(v)$ whenever box $u$ is weakly northwest of $v$. In the case~$k=1$ we call such a tableau a \dfn{standard barely set-valued tableau}. We use~$\SYTo(\lambda)$ to denote the set of standard barely set-valued tableaux of shape $\lambda$. Note that for a tableau $T \in \SYTo(\lambda)$, all but one of the boxes of $T$ have single numbers as their entries, and that special box has two numbers as its entry.

Set-valued tableaux were first introduced by Buch~\cite{buch2002littlewood} in his study of the $K$-theory of Grassmannians; in particular, they enter into the combinatorial definition of \dfn{(stable) Grothendieck polynomials}, which are $K$-theoretic extensions of \dfn{Schur functions}. Recently, set-valued tableaux have also appeared in a different geometric context: the algebraic geometry of curves. Let us review this recent algebro-geometric work. 

Let $C$ be a generic smooth curve of genus $g$. For parameters $r$ and $d$, the \dfn{Brill--Noether space} $G^r_d(C)$ is the moduli space of maps $C\to \mathbb{P}^{r}$ from $C$ into $r$-dimensional projective space of degree $d$. Brill--Noether theory~\cite{harris2009brill} is the study of this space $G^r_d(C)$. Define the number $\rho=\rho(g,r,d)$ by
\[ \rho \coloneqq g - (r+1)(g-d+r).\]
The Brill--Noether Theorem says that $G^r_d(C)$ is non-empty if and only if $\rho \geq 0$, and in this case $\rho$ is the dimension of~$G^r_d(C)$. 

We might be interested in finer numerical information about $G^r_d(C)$ than just its dimension. For example, when $\rho=0$, $G^r_d(C)$ is a $0$-dimensional variety, i.e., a collection of points, and the exact number of points is known to be
\begin{equation} \label{eqn:bn1}
g! \cdot \prod_{i=0}^{r} \frac{i!}{(g-d+r+i)!}.
\end{equation}
When $\rho=1$, $G^r_d(C)$ is itself a smooth curve, and the genus of this curve is known to be
\begin{equation} \label{eqn:bn2}
1+ \frac{(r+1)(g-d+r)}{g-d+2r+1} \cdot g! \cdot \prod_{i=0}^{r} \frac{i!}{(g-d+r+i)!}.
\end{equation}

It is easy to see from the hook length formula~\eqref{eqn:hlf} that~\eqref{eqn:bn1} is also equal to the number $\#\SYT(\rect{(r+1)}{(g-d+r)})$ of SYTs of $\rect{(r+1)}{(g-d+r)}$ \dfn{rectangular shape}. In some recent combinatorial approaches to Brill--Noether theory~\cite{cools2012tropical} the connection to tableaux is made explicit. Moreover, Chan et al.~\cite{chan2018genera} gave a tableau interpretation of the genus of $G^r_d(C)$ when $\rho=1$: they showed it is $1+\#\SYTo(\rect{(r+1)}{(g-d+r)})$.\footnote{Technically, Chan et al.~\cite{chan2018genera} did not use the language of set-valued tableaux. Indeed, the term ``barely set-valued tableau'' was introduced in the later paper of Reiner, Tenner, and Yong~\cite{reiner2018poset}. In~\cite{chan2018genera}, the genus of the Brill--Noether locus was related to the number of edges in a certain combinatorially defined graph (the ``Brill--Nother graph''). But it is easy to see that the edges of this graph correspond to barely set-valued tableaux: this point is clarified in \cite[Remark 6.5]{hopkins2017CDE}.} They were also able to give a new proof of~\eqref{eqn:bn2} by proving combinatorially that
\begin{equation} \label{eqn:bsv_rect}
 \#\SYTo(\rect{a}{b}) = \frac{ab}{a+b} \cdot (ab+1) \cdot \# \SYT(\rect{a}{b}),
\end{equation}
for any $a,b \geq 1$. In follow up work, Chan and Pflueger~\cite{chan2021euler} showed more generally that for any $\rho \geq 0$ the algebraic Euler characteristic of $G^r_d(C)$ is $(-1)^{\rho}$ times the number of standard set-valued tableaux of shape $\rect{(r+1)}{(g-d+r)}$ with entries $1,2,\ldots,g$. However, it is unclear whether there are product formulas for this quantity for all~$\rho \geq 0$, although there are known determinantal formulas~\cite{anderson2017kclasses, chan2021combinatorial}.

Independently of this work in Brill--Noether theory, Reiner, Tenner, and Yong~\cite{reiner2018poset} also investigated barely set-valued tableaux from the more traditional vantage point of symmetric functions. They proved a number of results extending~\eqref{eqn:bsv_rect}, and made a number of conjectures, and subsequently there has been a reasonable amount of research devoted to counting classes of barely set-valued tableaux~\cite{hopkins2017CDE, fan2019proof, kim2020enumeration}. We note that while~\eqref{eqn:hlf} applies to any partition, product formulas like~\eqref{eqn:bsv_rect} for $\#\SYTo(\lambda)$ exist only for very special shapes $\lambda$ like rectangles. Nevertheless, until all of this recent work there was no real reason to expect product formulas enumerating any classes of set-valued tableaux (see for example~\cite[Problem 3]{monical2019reduced}).

\subsection{Refinement over comaj for SYT}
In this paper we prove a $q$-analog of~\eqref{eqn:bsv_rect}. 

Before explaining our result, let us first recall the well-known $q$-analog of~\eqref{eqn:hlf}. Let $\lambda$ be a partition of $n$. For an SYT $T \in \SYT(\lambda)$, a \dfn{(natural) descent} of $T$ is an entry $i$ such that $i+1$ appears in a higher row than $i$.\footnote{The more common definition of a descent for an SYT, especially in the context of symmetric function theory, is an entry $i$ such that $i+1$ appears in a \emph{lower} row than $i$. Natural descents are descents with respect to a \dfn{natural labeling} of the poset of $\lambda$; the more common descents are descents with respect to a \dfn{Schur labeling}. The comajor index generating function for SYT  using the more common notion of descent (which also equals the corresponding \dfn{major index} generating function) is almost the same as~\eqref{eqn:qhlf} except that the term of $q^{-b(\lambda)}$ is not present; see~\cite[Proposition~7.19.11]{stanley1999ec2}.} We use $\Des(T)$ to denote the set of (natural) descents of~$T$. We define the \dfn{(natural) comajor index} of~$T$ to be $\comaj(T) \coloneqq \sum_{i \in \Des(T)} (n-i)$. The basic theory of $P$-partitions~\cite[\S3.15]{stanley1996ec1} (see also~\cite[Proof of Theorem 7.22.1]{stanley1999ec2}), and Stanley's \dfn{hook content formula}~\cite[Theorem 7.21.2]{stanley1999ec2}, then imply that
\begin{align} \label{eqn:qhlf}
\sum_{T \in \SYT(\lambda)} q^{\comaj(T)} &= (1-q) (1-q^2) \cdots (1-q^n) \cdot q^{-b(\lambda)} \cdot s_\lambda(1,q,q^2,\ldots) \\
\notag &= [n]_q! \cdot \prod_{u \in \lambda} \frac{1}{[h(u)]_q}.
\end{align}
Here $b(\lambda) \coloneqq\sum_{i=0}^{\infty} (i-1)\cdot \lambda_i$, $s_{\lambda}$ is the Schur function associated to $\lambda$, and we use the \dfn{$q$-number} notation $[k]_q \coloneqq\frac{(1-q^k)}{(1-q)} = 1+q+\cdots+q^{k-1}$ and $[n]_q! \coloneqq[n]_q[n-1]_q\cdots[1]_q$.

Now let us give the appropriate extensions of these notions to barely set-valued tableaux. Let $S \in \SYTo(\lambda)$. There is a unique box in $S$, call it $\ustar(S)$, containing two numbers. Set $\istar(S) \coloneqq\max(S(\ustar(S)))$. A \dfn{barely set-valued descent} of $S$ is then a number $i=1,2,\ldots,n+1$ such that $i+1$ appears in a higher row than $i$, \emph{except that}:
\begin{itemize}
\item $\istar(S)-1$ is \emph{never} a descent;
\item $\istar(S)$ is \emph{always} a descent.
\end{itemize}
We use $\Des^{+1}(S)$ to denote the set of descents of $S$. We define the \dfn{barely set-valued comajor index} of $S$ to be $\comaj^{+1}(S) \coloneqq \sum_{i \in \Des^{+1}(S)} (n+1-i)$.

Our first main result is the following $q$-analog of~\eqref{eqn:bsv_rect} (see \cref{cor:bsv_lin}):

\begin{thm} \label{thm:syt_intro}
For any $a,b\geq 1$,
\begin{align*}
\sum_{S \in \SYTo(\rect{a}{b})} q^{\comaj^{+1}(S)} &=\frac{[a]_q[b]_q}{[a+b]_q} \cdot [ab+1]_q \cdot \sum_{T\in \SYT(\rect{a}{b})} q^{\comaj(T)} \\
&=\frac{[a]_q[b]_q}{[a+b]_q} \cdot [ab+1]_q! \cdot \prod_{i=0}^{a-1} \frac{[i]_q!}{[b+i]_q!}.
\end{align*}
\end{thm}

\begin{example} \label{ex:syt}
Consider the case $a=b=2$. There are two elements of $\SYT(\rect{2}{2})$, which we show together with their descent sets and comajor indices in~\cref{tab:2x2_syt}. The generating function of their comajor indices is
\[ \sum_{T\in\SYT(\rect{2}{2})} q^{\comaj(T)} = q^2 + 1 = \frac{[4]_q!}{[3]_q[2]_q[2]_q[1]_q},\]
in agreement with~\eqref{eqn:qhlf}. Meanwhile, there are $10$ elements of $\SYTo(\rect{2}{2})$, which together with their descent sets and comajor indices are shown in~\cref{tab:2x2_bsv_syt}. The generating function of their comajor indices is
\[ \sum_{S\in\SYTo(\rect{2}{2})} q^{\comaj^{+1}(S)} =  q^5 + 2 q^4 + 2 q^3 + 2 q^2 + 2 q + 1 = \frac{[2]_q[2]_q}{[4]_q} \cdot [5]_q \cdot (q^2 + 1),\]
in agreement with \cref{thm:syt_intro}.
\end{example}
\begin{table}[ht]
\renewcommand{\arraystretch}{1.5}
\begin{tabular} {c | c | c}
$T$ & \parbox{0.5in}{\begin{center}\begin{ytableau} 1 & 2 \\ 3 & 4 \end{ytableau}\end{center}} & \parbox{0.5in}{\begin{center}\begin{ytableau} 1 & 3 \\ 2 & 4 \end{ytableau}\end{center}} \\ \hline
$\Des(T)$ & $\varnothing$ & $\{2\}$ \\ \hline
$\comaj(T)$ & $0$ & $2$
\end{tabular} 
\medskip
\caption{The tableaux in $\SYT(\rect{2}{2})$ with (natural) descents and comajor indices.} \label{tab:2x2_syt}
\end{table}

\begin{table}[ht]
\renewcommand{\arraystretch}{1.5}
\begin{tabular} {c | c | c | c | c | c}
$S$ & \parbox{0.5in}{\begin{center}\begin{ytableau} 1 & 2 \\ 3 & 4,5 \end{ytableau}\end{center}} & \parbox{0.5in}{\begin{center}\begin{ytableau} 1 & 3 \\ 2 & 4,5 \end{ytableau}\end{center}} & \parbox{0.5in}{\begin{center}\begin{ytableau} 1 & 2 \\ 3,4 & 5 \end{ytableau}\end{center}} & \parbox{0.5in}{\begin{center}\begin{ytableau} 1 & 3 \\ 2,4 & 5 \end{ytableau}\end{center}} & \parbox{0.5in}{\begin{center}\begin{ytableau} 1 & 4 \\ 2,3 & 5 \end{ytableau}\end{center}} \\ \hline
$\Des^{+1}(S)$ & $\{5\}$ & $\{2,5\}$ & $\{4\}$ & $\{2,4\}$ & $\{3\}$ \\ \hline
$\comaj^{+1}(S)$ & $0$ & $3$ & $1$ & $4$ & $2$
\end{tabular} 
\medskip \hrule \medskip
\begin{tabular} {c | c | c | c | c | c}
$S$ & \parbox{0.5in}{\begin{center}\begin{ytableau} 1 & 2,3 \\ 4 & 5 \end{ytableau}\end{center}} & \parbox{0.5in}{\begin{center}\begin{ytableau} 1 & 2,4 \\ 3 & 5 \end{ytableau}\end{center}} & \parbox{0.5in}{\begin{center}\begin{ytableau} 1 & 3,4 \\ 2 & 5 \end{ytableau}\end{center}} & \parbox{0.5in}{\begin{center}\begin{ytableau} 1,2 & 3 \\ 4 & 5 \end{ytableau}\end{center}} & \parbox{0.5in}{\begin{center}\begin{ytableau} 1,2 & 4 \\ 3 & 5 \end{ytableau}\end{center}} \\ \hline
$\Des^{+1}(S)$ & $\{3\}$ & $\{4\}$ & $\{2,4\}$ & $\{2\}$ & $\{2,3\}$ \\ \hline
$\comaj^{+1}(S)$ & $2$ & $1$ & $4$ & $3$ & $5$
\end{tabular} 
\medskip
\caption{The tableaux in $\SYTo(\rect{2}{2})$ with barely set-valued descents and comajor indices.} \label{tab:2x2_bsv_syt}
\end{table}

\begin{remark} \label{rem:rows}
For $S \in \SYTo(\lambda)$, let $\rstar(S)$ denote the row that $\ustar(S)$ occurs in. Extending~\cref{thm:syt_intro}, we will in fact prove (see Corollary \ref{cor:bsv_lin_refined})
\[ \sum_{S \in \SYTo(\rect{a}{b})} q^{\comaj^{+1}(S)}t^{\rstar(S)-1} = \frac{[a]_{q,t} [b]_q}{[a+b]_q} \cdot [ab+1]_q \cdot \sum_{T\in \SYT(\rect{a}{b})}q^{\comaj(T)},\]
where $[k]_{q,t} \coloneqq \frac{(t^k-q^k)}{(t-q)} = (t^{k-1}+qt^{k-2}+\cdots+q^{k-2}t+q^{k-1})$. We encourage the reader to verify this for \cref{ex:syt}.
\end{remark}

\subsection{Barely set-valued plane partitions} 
Recall that an $\rect{a}{b}$ \dfn{plane partition} is an $\rect{a}{b}$ array of nonnegative integers that is weakly decreasing in rows and columns. For a plane partition $\pi$, the \dfn{size} $|\pi|$ of~$\pi$ is the sum of all the entries of~$\pi$. Let us use $\PP_m(\rect{a}{b})$ to denote the set of $\rect{a}{b}$ plane partitions with entries in~$\{0,1,\ldots,m\}$. MacMahon~\cite{macmahon1915combinatory} famously proved the following product formula for the size generating function of such plane partitions:
\begin{equation} \label{eqn:macmahon}
\sum_{\pi \in \PP_m(\rect{a}{b})} q^{|\pi|} = \prod_{i=1}^{a}\prod_{j=1}^{b} \frac{[i+j+m-1]_q}{[i+j-1]_q}.
\end{equation}

We define a \dfn{barely set-valued plane partition} to be like a plane partition except that exactly one box in the array contains two different numbers, while all the other boxes have a single number. Let~$\PPo_m(\rect{a}{b})$ denote the set of barely set-valued $\rect{a}{b}$ plane partitions with entries in $\{0,1,\ldots,m\}$. For a barely set-valued plane partition~$\tau$ we again use~$|\tau|$ to denote the sum of the entries of $\tau$. We prove the following barely set-valued analog of MacMahon's formula~\eqref{eqn:macmahon} (see \cref{cor:bsv_rpp}):

\begin{thm} \label{thm:pp_intro}
For any $a,b,m\geq 1$,
\[\sum_{\tau \in \PPo_m(\rect{a}{b})} \hspace{-0.4cm} q^{|\tau|-1} = \frac{[a]_q[b]_q}{[a+b]_q} \cdot [m]_q \cdot \hspace{-0.4cm} \sum_{\pi \in \PP_m(\rect{a}{b})} \hspace{-0.4cm} q^{|\pi|} = \frac{[a]_q[b]_q}{[a+b]_q} \cdot [m]_q \cdot \prod_{i=1}^{a}\prod_{j=1}^{b} \frac{[i+j+m-1]_q}{[i+j-1]_q}.\]
\end{thm}

\begin{example}\label{ex:pp}
Consider the case $a=b=2$ and $m=1$. There are six elements of $\PP_1(\rect{2}{2})$, which we show together with their sizes in~\cref{tab:2x2_pp}. The generating function of their sizes is
\[ \sum_{\pi \in\PP_1(\rect{2}{2})} q^{|\pi|} = q^4 + q^3 + 2q^2 + q + 1 = \frac{[4]_q[3]_q[3]_q[2]_q}{[3]_q[2]_q[2]_q[1]_q},\]
in agreement with~\eqref{eqn:macmahon}. Meanwhile, there are also six elements of $\PPo_1(\rect{2}{2})$, which together with their sizes are shown in~\cref{tab:2x2_bsv_pp}. The generating function of their sizes (minus one) is
\[ \sum_{\tau \in\PP_1^{+1}(\rect{2}{2})} q^{|\tau|-1} =  q^3+2q^2+2q+1 = \frac{[2]_q[2]_q}{[4]_q} \cdot [1]_q \cdot (q^4 + q^3 + 2q^2 + q + 1),\]
in agreement with \cref{thm:pp_intro}.
\end{example}

\begin{table}[ht]
\renewcommand{\arraystretch}{1.5}
\begin{tabular} {c | c | c | c | c | c | c}
$\pi$ & \parbox{0.5in}{\begin{center}\begin{ytableau} 0 & 0 \\ 0 & 0 \end{ytableau}\end{center}} & \parbox{0.5in}{\begin{center}\begin{ytableau} 1 & 0 \\ 0 & 0 \end{ytableau}\end{center}} & \parbox{0.5in}{\begin{center}\begin{ytableau} 1 & 1 \\ 0 & 0 \end{ytableau}\end{center}} & \parbox{0.5in}{\begin{center}\begin{ytableau} 1 & 0 \\ 1 & 0 \end{ytableau}\end{center}} & \parbox{0.5in}{\begin{center}\begin{ytableau} 1 & 1 \\ 1 & 0 \end{ytableau}\end{center}} & \parbox{0.5in}{\begin{center}\begin{ytableau} 1 & 1 \\ 1 & 1 \end{ytableau}\end{center}} \\ \hline
$|\pi|$ & $0$  & $1$ & $2$ & $2$ & $3$ & $4$ \\
\end{tabular} 
\medskip
\caption{The plane partitions in $\PP_1(\rect{2}{2})$ with sizes.} \label{tab:2x2_pp}
\end{table}

\begin{table}[ht]
\renewcommand{\arraystretch}{1.5}
\begin{tabular} {c | c | c | c | c | c | c}
$\tau$ & \parbox{0.5in}{\begin{center}\begin{ytableau} 1,0 & 0 \\ 0 & 0 \end{ytableau}\end{center}} & \parbox{0.5in}{\begin{center}\begin{ytableau} 1 & 1,0 \\ 0 & 0 \end{ytableau}\end{center}} & \parbox{0.5in}{\begin{center}\begin{ytableau} 1 & 0 \\ 1,0 & 0 \end{ytableau}\end{center}} & \parbox{0.5in}{\begin{center}\begin{ytableau} 1 & 1,0 \\ 1 & 0 \end{ytableau}\end{center}} & \parbox{0.5in}{\begin{center}\begin{ytableau} 1 & 1 \\ 1,0 & 0 \end{ytableau}\end{center}} & \parbox{0.5in}{\begin{center}\begin{ytableau} 1 & 1 \\ 1 & 1,0 \end{ytableau}\end{center}} \\ \hline
$|\tau|$ & $1$  & $2$ & $2$ & $3$ & $3$ & $4$ \\
\end{tabular} 
\medskip
\caption{The plane partitions in $\PPo_1(\rect{2}{2})$ with sizes.} \label{tab:2x2_bsv_pp}
\end{table}

There is a precise sense in which, for $\lambda=\rect{a}{b}$, the formula in~\eqref{eqn:qhlf} is the $m\to\infty$ limit of \eqref{eqn:macmahon}. Similarly, \cref{thm:syt_intro} is the the $m\to\infty$ limit of \cref{thm:pp_intro}.

\subsection{Outline of proof}
Our proofs of \cref{thm:syt_intro,thm:pp_intro} follow an approach of Reiner--Tenner--Yong~\cite{reiner2018poset} and Chan et al.~\cite{chan2017expected}. The first step is to reformulate everything we have discussed so far (tableaux, plane partitions, et cetera) in the more general language of posets. Indeed, even though ultimately we obtain product formulas only for very special shapes like the rectangles, many of the results we prove along the way hold for arbitrary posets $P$.

The meat of our proof consists of defining and studying certain probability distributions $\murpp{m}{q}$ and $\mulin{q}$ (depending on a real parameter $q>0$) on the set $\J(P)$ of \dfn{order ideals} of a finite poset $P$. As mentioned, this general strategy for counting barely set-valued fillings using probability distributions was essentially carried out in~\cite{chan2017expected,reiner2018poset}. But those papers only dealt with the case $q=1$. Our main contribution here is to introduce the $q$ parameter. 

The distributions we consider have two important properties. The first important property is that they are related to barely set-valued fillings. More precisely, we can $q$-count barely set-valued fillings of $P$ by computing the expectation of the \dfn{down-degree} statistic $\ddeg\colon \J(P)\to\R$ under these distributions. For example, in the case when $P=\rect{a}{b}$ is the rectangle poset, we have
\begin{align} \label{eqn:ex_intro}
\E_{\murpp{m}{q}}(\ddeg) &= \frac{\sum_{\tau \in \PPo_m(\rect{a}{b})} q^{|\tau|-1}}{[m]_q \cdot \sum_{\pi \in \PP_m(\rect{a}{b})} q^{|\pi|}}, \\
\notag \E_{\mulin{q}}(\ddeg) &= \frac{\sum_{S \in \SYTo(\rect{a}{b})} q^{\comaj^{+1}(S)}}{[ab+1]_q \cdot \sum_{T\in \SYT(\rect{a}{b})} q^{\comaj(T)}}.
\end{align}
The second important property of these probability distributions is a certain symmetry property called \dfn{$q$-toggle symmetry}. 

As we will prove below, both of these two important properties of the probability distributions $\murpp{m}{q}$ and $\mulin{q}$ hold for arbitrary posets $P$. But what makes the rectangle special is a third property: for any $q$-toggle symmetric probability distribution $\mu$ on $\J(\rect{a}{b})$, we have
\[\E_{\mu}(\ddeg) = \frac{[a]_q[b]_q}{[a+b]_q}.\]
Combining all of the above yields the elegant formulas in \cref{thm:syt_intro,thm:pp_intro}. 

To establish the third property for the rectangle (that all $q$-toggle symmetric distributions give the same expected down-degree), we rely crucially on the very recent paper~\cite{defant2021homomesy}. Actually, this third property holds for a few other special posets, including for the \dfn{shifted staircase} and, more generally, for all the \dfn{minuscule posets}. Hence, as we explain, our product formulas for barely set-valued fillings extend to all minuscule posets.

\smallskip

The rest of the paper is structured as follows. In \cref{sec:dists} we review some basic notions related to posets including order ideals, linear extensions, $P$-partitions, et cetera, and explain how the tableaux constructions discussed in this introduction are special cases of these more general notions. We then define several probability distributions on the set of order ideals of a poset. These distributions all depend on a parameter $q > 0$. In \cref{sec:tog_sym} we single out the class of ``$q$-toggle-symmetric'' probability distributions on the order ideals of a poset. These are the distributions where we are~$q$ times as likely to be able to toggle in as toggle out every element of the poset. We show that all the distributions from~\cref{sec:dists} are $q$-toggle-symmetric. In~\cref{sec:proof_rect} we connect the distributions to barely set-valued fillings. We explain how such fillings can be counted (and $q$-counted) using expectations of random variables, as in \eqref{eqn:ex_intro}. Finally, we compute these expectations in the case of the rectangle by appealing to the $q$-toggle-symmetry of the relevant distributions, as well as a crucial result proved in~\cite{defant2021homomesy}. In~\cref{sec:sstair} we extend these enumerative results to the shifted staircase and all the minuscule posets. In~\cref{sec:final} we end with some final remarks and open questions.

\subsection*{Acknowledgements}
The first author was supported by \emph{National Science Foundation} grant \#1802920. The second and third author were funded by \emph{Swedish Research Council} grant 2018-05218. We thank the organizers of the 2020 Banff International Research Stations (BIRS) online workshop on Dynamical Algebraic Combinatorics for giving us a chance to collaborate. In particular, we thank Jim Propp and the other authors of~\cite{defant2021homomesy}, because this project ultimately grew out of the investigation of ``$q$-rowmotion,'' which he initiated and which is carried out in that paper. Finally, we thank the anonymous referees for their careful reading of our paper and useful comments which improved the exposition.

\section{Distributions on order ideals} \label{sec:dists}

In this section we define several probability distributions on the set of order ideals of a finite poset. The connection to the discussion of our results in~\cref{sec:intro} may not be totally clear, and may not become clear until \cref{sec:proof_rect}. We hope the reader bears with us.

\subsection{Background on posets}

We assume the reader is familiar with the basic notions and notations from the theory of posets as laid out for instance in~\cite[Chapter~3]{stanley1996ec1}. In particular we use the notation $[n] \coloneqq \{1,2,\ldots,n\}$, which we also view as a chain poset with the usual total order.
Throughout the rest of the paper, $P$ will be a poset on $n$ elements.
An \dfn{order ideal} of $P$ is a downwards-closed subset: i.e. a subset $I\subseteq P$ with $y\in I$ and $x \leq y \in P$ implies $x \in I$. The set of order ideals of $P$ is denoted $\J(P)$. Although we will not need this, we note that $\J(P)$ is a distributive lattice with respect to the inclusion order.

Let $\N\coloneqq\{0,1,\ldots,\}$ denote the set of natural numbers. A \dfn{reverse $P$-partition} is a map $\pi\colon P \to \N$ that is order-preserving: $p\leq p'$ implies $\pi(p)\leq \pi(p')$.\footnote{Traditionally \dfn{$P$-partitions}, which are order-reversing maps, are studied. We prefer to work with order-preserving maps. This means that in many places our convention is the opposite of that of~\cite{stanley1996ec1} (e.g., we always use comajor index rather than major index). Ultimately the translation is only a matter of replacing $P$ by its \dfn{dual poset} $P^*$.} We use $\RPP(P)$ to denote the set of reverse $P$-partitions; and we use $\RPP_m(P)$ to denote the set of reverse $P$-partitions with $\pi(P)\subseteq \{0,1,\ldots,m\}$. Note that we have a canonical bijection $\RPP_1(P) \xrightarrow{\sim} \J(P)$ where $\pi \mapsto \pi^{-1}(0)$. A related fact that we will use constantly is that for any $\pi \in \RPP(P)$ and $i \in \N$, we have $\pi^{-1}(\{0,1,\ldots,i\}) \in \J(P)$.

A \dfn{linear extension} of $P$ is a bijection $T\colon P \to [n]$ that is order-preserving. We use~$\Lin(P)$ to denote the set of linear extensions of $P$.
From now on we assume that~$P$ comes with a fixed ``reference'' linear extension $\omega$; such an $\omega$ is called a \dfn{natural labeling} of $P$.

We use $\mathfrak{S}_n$ to denote the symmetric group of permutations of $[n]$. Via $\omega$, every linear extension of $P$ is identified with a permutation in $\mathfrak{S}_n$: $T \in \Lin(P)$ is identified with the permutation $w_T$ whose one-line notation is 
\[\omega(T^{-1}(1)),\omega(T^{-1}(2)),\ldots,\omega(T^{-1}(n)).\]

Recall that for a permutation $w=w_1\cdots w_n \in \mathfrak{S}_n$, a \dfn{descent} of $w$ is an index $i=1,2,\ldots,n-1$ for which $w_i > w_{i+1}$. We thus define a \dfn{descent} of a linear extension $T \in \Lin(P)$ to be a descent of $w_T$. We use $\Des(T)$ to denote the set of descents of $T$.

\smallskip

Now let us explain how these poset-theoretic concepts relate to the results on tableaux we discussed in~\cref{sec:intro}. We assume that the reader is familiar with the basics regarding integer partitions, Young diagrams, et cetera, as laid out for instance in~\cite[Chapter~1]{stanley1996ec1}. We use ``English notation'' for Young diagrams, meaning that, for example, the partition $\lambda = (4,3,3)$ has Young diagram
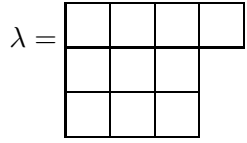
\[ \lambda = \ydiagram{4,3,3}\]
We use matrix coordinates for the boxes of a Young diagram, so that the northwest-most box is $(1,1)$, the box immediately to its east is $(1,2)$, and so on.

Any partition $\lambda$ of $n$ determines a poset $\shape{\lambda}$ whose elements are the boxes of the Young diagram of $\lambda$, and with $u\leq u'$ if $u$ is weakly northwest of $u'$. We often refer to posets of the form $\shape{\lambda}$ as \dfn{shapes}. Observe that SYTs of shape $\lambda$ are the same as linear extensions of $\shape{\lambda}$. We always assume that $\shape{\lambda}$ is given the natural labeling~$\omega$ where the 1st row of $\lambda$ gets the values $1,2,\ldots,\lambda_1$, the 2nd row $\lambda_1+1,\ldots,\lambda_2$, and so on. With this choice of $\omega$, (natural) descents of SYTs as we defined them in~\cref{sec:intro} are the same as descents of linear extensions of $\shape{\lambda}$ as we just defined them.

By abuse of notation, for $a,b\geq 1$ we use $\rect{a}{b}$ to denote both the rectangular partition $\rect{a}{b} \coloneqq (b^a)$ and the associated \dfn{rectangle poset} $\rect{a}{b} \coloneqq \shape{\rect{a}{b}}$. Note that another notation we could use for the rectangle poset is $[a]\times[b]$, i.e., the (Cartesian) product of two chains.

\smallskip

We return to discussing concepts for general posets $P$ on $n$ elements. For any permutation $w = w_1w_2\cdots w_n \in \mathfrak{S}_n$, a function $f: P \to \N$ is said to be \dfn{$w$-compatible} if
\begin{enumerate}
    \item $f(\omega^{-1}(w_1)) \leq f(\omega^{-1}(w_2)) \leq \cdots \leq f(\omega^{-1}(w_n))$
    \item $f(\omega^{-1}(w_i)) < f(\omega^{-1}(w_{i+1}))$ if $w_i > w_{i+1}$.
\end{enumerate}
The following is the basic result that makes the theory of $P$-partitions work (see~\cite[Theorem~3.15.7]{stanley1996ec1}):

\begin{lemma} \label{lem:w_compatible}
A function $f: P \to \N$ is a reverse $P$-partition if and only if $f$ is compatible with $w_T$ for some (unique) linear extension $T \in \Lin(P)$.
\end{lemma}

There are many important consequences of~\cref{lem:w_compatible}. One consequence is that the multiset of descent sets $\Des(T)$ over all $T \in \Lin(P)$ does not depend on the particular natural labeling $\omega$ we have fixed. Another consequence is the following. For a reverse $P$-partition $\pi \in \RPP(P)$, define its \dfn{size} to be $|\pi| \coloneqq \sum_{p\in P} \pi(p)$. For a linear extension $T \in \Lin(P)$, define its \dfn{comajor index} to be $\comaj(T) \coloneqq \sum_{i\in \Des(T)} (n-i)$. Then we have (see~\cite[Lemma 3.15.3]{stanley1996ec1}):
\begin{equation} \label{eqn:rpp_gf}
\sum_{\pi \in \RPP(P)} q^{|\pi|} = \frac{\sum_{T \in \Lin(P)}q^{\comaj(T)}}{(1-q)(1-q^2)\cdots(1-q^{n})}
\end{equation}
as an equality of formal power series in the parameter $q$.

\subsection{The \texorpdfstring{$q$}{q}-uniform distribution}

We will now define several probability distributions on $\J(P)$. These probability distributions will all depend on a real number parameter $q > 0$. In the case when $q=1$, they were all considered in~\cite{chan2017expected}.
So, from now on fix such a real number parameter $q > 0$. 

By a probability distribution $\mu$ on $\J(P)$ we just mean a function $\mu\colon\J(P) \to \R$ satisfying the usual axioms:
\begin{itemize}
\item $\mu(I)\geq 0$ for all $I \in \J(P)$;
\item $\sum_{I \in \J(P)}\mu(I)=1$.
\end{itemize}
We often say simply ``distribution'' instead of ``probability distribution.''
The first distribution we consider is $\muuni{q}\colon\J(P) \to \R$, given by
\begin{equation} \label{eqn:dist_uni}
\muuni{q}(I) \coloneqq\frac{q^{\#(P\setminus I)}}{\sum_{I\in \J(P)}q^{\#(P\setminus I)}}
\end{equation}
for all $I\in \J(P)$. More intuitively: we select the order ideal $I$ with probability proportional to~$q^{\# (P \setminus I)}$. Evidently $\muuni{q}$ becomes the uniform distribution in the case $q=1$, and hence the name.

\subsection{The reverse \texorpdfstring{$P$}{P}-partition distributions}

Now we will define a family of distributions in terms of reverse $P$-partitions. 
Let $m \geq 1$ be an integer. Define the distribution $\murpp{q}{m}\colon\J(P) \to \R$ by
\begin{equation}  \label{eqn:dist_rpp}
\murpp{q}{m}(I) \coloneqq\frac{\displaystyle \sum_{\pi \in \RPP_m(P)} \sum_{k=0}^{m-1} q^{|\pi|+k} \cdot \bracket{\pi^{-1}(\{0,1,\ldots,k\}) = I}}{\displaystyle [m]_q\cdot\sum_{\pi \in \RPP_m(P)} q^{|\pi|}}
\end{equation}
for all $I\in\J(P)$, where for a proposition $\phi$ we use the notation
\[\bracket{\phi} \coloneqq\begin{cases} 1 &\textrm{if $\phi$ is true}; \\ 0 &\textrm{otherwise}.\end{cases}\] 
More intuitively, $\murpp{q}{m}(I)$ is the probability that we obtain $I$ from the following procedure:
we select $\pi \in \RPP_m(P)$ with probability proportional to $q^{|\pi|}$; we select $k \in \{0,1,\ldots,m-1\}$ independently with probability proportional to $q^{k}$; and then we select the order ideal $\pi^{-1}(\{0,1,\ldots,k\})$.
It is clear that $\murpp{q}{m}$ really defines a distribution, i.e., that $\sum_{I\in \J(P)}\murpp{q}{m}(I)=1$.
It is also clear that $\murpp{q}{1}=\muuni{q}$, i.e., that \eqref{eqn:dist_rpp} reduces to \eqref{eqn:dist_uni}  in the case $m=1$.

\begin{example}\label{ex:rpp_probs}
In Table \ref{tab:2x2_pp} we listed the six elements of $\PP_1(\rect{2}{2})$ along with their sizes. Notice that, after rotating by $180^{\circ}$, these elements are also the six elements of $\RPP_1(\rect{2}{2})$. By Example \ref{ex:pp}, we have
\[\sum_{\pi \in \RPP_1(\rect{2}{2})}q^{|\pi|} = q^4 + q^3 + 2q^2 + q + 1.\]
We can now compute $\murpp{q}{1}(I)$ for all six order ideals of $\rect{2}{2}$: 
\[\begin{minipage}{0.49\textwidth}
\ytableausetup{smalltableaux}
\begin{align*}
\murpp{q}{1}\left( \varnothing \right) &= \frac{q^4}{q^{4}+q^3+2q^2+q+1}\\
\murpp{q}{1}\left( \parbox{0.4cm}{\centering\ydiagram{1}} \right) &= \frac{q^3}{q^{4}+q^3+2q^2+q+1}\\
\murpp{q}{1}\left( \parbox{0.8cm}{\centering\ydiagram{2}} \right) &= \frac{q^2}{q^{4}+q^3+2q^2+q+1}
\end{align*}
\end{minipage}
\begin{minipage}{0.49\textwidth}
\begin{align*}
\murpp{q}{1}\left( \parbox{0.4cm}{\centering\ydiagram{1,1}} \right) &= \frac{q^2}{q^{4}+q^3+2q^2+q+1}\\
\murpp{q}{1}\left(  \parbox{0.8cm}{\centering\ydiagram{2,1}} \right) &= \frac{q}{q^{4}+q^3+2q^2+q+1}\\
\murpp{q}{1}\left(  \parbox{0.8cm}{\centering\ydiagram{2,2}} \right) &= \frac{1}{q^{4}+q^3+2q^2+q+1}
\end{align*}
\ytableausetup{nosmalltableaux}
\end{minipage}\]
Notice that these probabilities agree with the $q$-uniform distribution, as we just observed.
\end{example}

For convenience, we introduce the following notation that we will use frequently in the next several sections. For any $T \in \Lin(P)$ and $i=0,1,\ldots,n$, we define
\[\comaj(T,i) \coloneqq \displaystyle\sum_{j \in \Des(T)\cup \{i\}}n-j.\]
Also, we recall that the \dfn{$q$-binomial coefficient} $\qbinom{n}{k}_q$ is $\qbinom{n}{k}_q \coloneqq\frac{[n]_q!}{[k]_q![n-k]_q!}$ if $n \geq k$, and zero otherwise.

The sums in the numerator of~\eqref{eqn:dist_rpp} are over sets whose cardinality is growing as~$m$ grows. It would be preferable, especially for the analysis of the limit as $m \to \infty$, to be able to express this numerator as a sum over a fixed set, independent of $m$. In fact, we can do that, as the next theorem explains.

\begin{thm}\label{thm:m-weight}
For any poset $P$ and any $I \in \J(P)$,
\begin{gather} \label{eqn:m-weight}
\sum_{\pi \in \RPP_m(P)} \sum_{k=0}^{m-1} q^{|\pi|+k} \cdot \bracket{\pi^{-1}(\{0,1,\ldots,k\}) = I} \\ 
\notag= \sum_{T\in\Lin(P)}\!\!\! \wt_m(T,\#I) \cdot \bracket{T^{-1}(\{1,\ldots,\#I\}) = I},
\end{gather}
where 
\[\displaystyle \wt_m(T,i) \coloneqq q^{\comaj(T,i)} \cdot q^{\#\{j \in \Des(T) \colon  j < i\}} \cdot \qbinom{m+n-\#(\Des(T)\setminus \{i\})}{n+1}_q,\] 
for all $T \in \Lin(P)$ and $i=0,1,\ldots,n$.
\end{thm}

\begin{proof}
Fix an $I \in \J(P)$. By Lemma \ref{lem:w_compatible}, we know that the set $\RPP_m(P)$ can be partitioned according to the $w_T$-compatibility of $\pi$ for $T \in \Lin$. Thus, 
the left-hand side of \eqref{eqn:m-weight} can then be rewritten as
\[\sum_{T \in \Lin(P)} \; \sum_{\substack{\pi \in \RPP_m(P)\\ \pi \textrm{ is $w_T$-compatible}}} \; \sum_{k=0}^{m-1} q^{|\pi|+k} \cdot \bracket{\pi^{-1}(\{0,1,\ldots,k\}) = I}\]
Let us focus on the contribution of the term corresponding to some fixed $T \in \Lin(P)$ in the above sum. By writing $w_T=w_1\cdots w_n$, and considering each $k\in \{0,\ldots,m-1\}$ one at a time, we have
\begin{equation}\label{eqn:qcountLHS}
\sum_{\substack{\pi \in \RPP_m(P),\\ \pi \textrm{ is $w_T$-compatible}}} \sum_{k=0}^{m-1} q^{|\pi|+k} \cdot \bracket{\pi^{-1}(\{0,1,\ldots,k\}) = I}=\sum_{\substack{0\leq \pi(w_1)\leq \cdots \leq \pi(w_{\#I})\leq k <\\ \pi(w_{\#I+1}) \leq \cdots \leq \pi(w_n) \leq m}}q^{|\pi|}q^{k},\end{equation}
where $\pi(w_j) < \pi(w_{j+1})$ whenever $w_j > w_{j+1}$ in the sum on the right-hand side. By the definition of reverse $P$-partitions, the sum is non-empty if and only if the first $\#I$ entries of $w$ are the elements of $I$. In such cases we can write $w_T$ as the concatenation $w_{T_1}\cdot w_{T_2}$ where $T_1$ is a linear extension of $I$ and $T_2$ is a linear extension of $P\setminus I$ using the alphabet $\{\#I+1,\dots,n\}$.

Our goal now is to ``compress'' the sum in the right-hand side of~\eqref{eqn:qcountLHS} by reindexing in order to remove all the strict inequality conditions. To that end, for each $j=1,\ldots,n$, let $x_j \coloneqq \pi(w_j)$ and define $\delta(j)$ to be the number of strict inequalities to the left of $x_j$ that are required for the conditions of $w$-compatibility, along with the strict inequality $k < x_{\#I +1}$. More precisely,
\[\delta(j) = \begin{cases} \#\{\ell \in \Des(T_1) \colon \ell < j\}, & 1\leq j \leq \#I\\
\#\{\ell \in \Des(T_2) \colon \ell < j\} + \#\Des(T_1) + 1, & \#I+1 \leq j \leq n.\end{cases}\]
Now we define a new set of variables and constants:
\begin{align*}x_j' &\coloneqq x_j - \delta(j)\\
k' &\coloneqq k - \#\Des(T_1)\\ 
m' &\coloneqq m - \#\Des(T_1) - \#\Des(T_2) - 1.\end{align*} 
We can then rewrite the index of summation in the right-hand side of \eqref{eqn:qcountLHS} as 
\[0\leq x_1'\leq x_2' \leq \cdots \leq x_{\#I}' \leq k' \leq x_{\#I+1}' \leq \cdots \leq x_n' \leq m'.\]
And observe that
\[q^{|\pi|} = q^{\sum x_j'}q^{\sum\delta(j)} \quad \text { and }\quad q^{k} = q^{k'}q^{\#\Des(T_1)},\]
so that
\begin{align*}
\sum_{\substack{0\leq x_1 \leq \cdots \leq x_{\#I} \leq k\\ \leq x_{\#I+1} \leq \cdots \leq x_n \leq m'}}q^{\sum x_j}q^{k}&= \left(\sum_{\substack{0\leq x_1' \leq \cdots \leq x_{\#I}' \leq k'\\ \leq x_{\#I+1}' \leq \cdots \leq x_n' \leq m'}}q^{\sum x_j'}q^{k'}\right)q^{\sum\delta(j)}q^{\#\Des(T_1)}\\
&= \left(\sum_{\pi \in \RPP_{m'}([n+1])}q^{|\pi|}\right)q^{\sum\delta(j)}q^{\#\Des(T_1)}\\
&= \qbinom{m'+n+1}{n+1}_q q^{\sum\delta(j)}q^{\#\Des(T_1)}.
\end{align*}
The last equality follows from the well-known combinatorial interpretation of the $q$-binomial coefficient $\qbinom{a+b}{b}_q$ as the generating function for partitions with at most $b$ parts and with each part less than or equal to $a$ (see, e.g.,~\cite[Proposition 1.7.3]{stanley1996ec1}). It can also be obtained from MacMahon's formula (\ref{eqn:macmahon}).

To finish the proof, we need to show that
\[ q^{\sum\delta(j)}q^{\#\Des(T_1)}\qbinom{m'+n+1}{n+1}_q = q^{\comaj(T,i)} q^{\#\{j \in \Des(T) \colon  j < i\}} \qbinom{m+n-\#(\Des(T)\setminus \{i\})}{n+1}_q,\]
where $i=\#I$. First, we need to compute $\displaystyle\sum_{j=1}^n\delta(j)$. So let $\alpha_1 <\dots < \alpha_{\ell_1}$ be the descents of $T_1$ and let $\beta_1 < \cdots < \beta_{\ell_2}$ be the descents of $T_2$. We have
\[\delta(j) = \begin{cases} 0, & j \leq \alpha_1\\ r, & \alpha_r < j \leq \alpha_{r+1},\; 1\leq r < \ell_1\\
\ell_1, & \alpha_{\ell_1} < j \leq \#I\\
\ell_1+1, & \#I<j\leq \beta_1\\
\ell_1 + r +1, & \beta_r < j \leq \beta_{r+1}, \; 1\leq r < \ell_2\\
\ell_1 + \ell_2 + 1, & \beta_{\ell_2} < j \leq n\end{cases}.\]
This means that we can write $\displaystyle\sum_{j=1}^n\delta(j)$ as:
\begin{align*}
&(\alpha_2 - \alpha_1) + 2(\alpha_3-\alpha_2) + \cdots + \ell_1(\#I - \alpha_{\ell_1})\\
&+ (\ell_1+1)(\beta_1 - \#I) + (\ell_1+2)(\beta_2 - \beta_1) + \cdots + (\ell_1+\ell_2 +1)(n - \beta_{\ell_2}).
\end{align*}
which simplifies to
$$(\ell_1 + \ell_2 +1)n - \sum_{r=1}^{\ell_1}\alpha_r - \sum_{r=1}^{\ell_2}\beta_r - \#I = (n-\#I) + \sum_{j \in \Des(T_1)}(n-j) + \sum_{j \in \Des(T_2)}(n-j).$$
Hence,
\begin{align*}q^{\sum\delta(j)}\cdot q^{\#\Des(T_1)}&= q^{(n-\#I)}\cdot q^{\#\Des(T_1)}\prod_{j \in \Des(T_1)}q^{n-j}\cdot\prod_{j \in \Des(T_2)}q^{n-j}\\
&= q^{(n-\#I)}\cdot q^{\#\Des(T_1)}\prod_{\substack{j \in \Des(T)\\ j <\#I}}q^{n-j}\cdot\prod_{\substack{j \in \Des(T)\\ j>\#I}}q^{n-j}\\
&= q^{\comaj(T,\#I)}\cdot q^{\#\{j \in \Des(T) \colon j < \#I\}},\end{align*}
since $j \in \Des(T_1)$ if and only if $j \in \Des(T)$ and $j < \#I$ and $j \in \Des(T_2)$ if and only if $j \in \Des(T)$ and $j > \#I$.
Moreover, we have
\[m' = m - \#\Des(T_1) - \#\Des(T_2) - 1 = m - \#(\Des(T)\setminus\{i\})-1.\]
Thus,
\[\qbinom{m'+n+1}{n+1}_q = \qbinom{m+n - \#(\Des(T)\setminus\{i\})}{n+1}_q.\]
This completes the proof.
\end{proof}

\begin{example}
Here we show an example of computing $\wt_m(T,i)$ when $P=\rect{3}{3}$. Let $i=3$, and let $T$ be the following SYT of shape $\rect{3}{3}$:
\begin{center}
\begin{ytableau}
*(yellow) 1 & *(yellow) 3 & 6 \\
*(yellow) 2 & 5 & 8\\
4 & 7 & 9 
\end{ytableau}
\end{center}
We have $D(T)=\{2,4,5,7\}$ and thus $\comaj(T,3)=7+6+5+4+2=24$, so
\[\wt_m(T,i) = q^{24}q^1\qbinom{m+9-4}{9+1}_q = q^{25}\qbinom{m+5}{10}_q.\]
Observe that $I \coloneqq T^{-1}(\{1,2,3\}) \in \J(P)$ is the staircase shape $(2,1)$ highlighted in yellow above. Hence, according to \cref{thm:m-weight}, we get a contribution of $q^{25}\qbinom{m+5}{10}_q$ to the weight of $I=(2,1)$ in $\murpp{q}{m}$ from this particular choice of $T$ and $i$.
\end{example}

\subsection{The linear extension distribution}

Now we define a distribution in terms of the linear extensions of $P$, and their descents. We will prove in Proposition \ref{prop:dist_lin} that this is indeed a distribution.

Define the distribution $\mulin{q} \colon \J(P)\to \R$ by
\begin{equation} \label{eqn:dist_lin}
\mulin{q}(I) \coloneqq\frac{\displaystyle \sum_{T \in \Lin(P)}\wt(T,\#I) \cdot \bracket{T^{-1}(\{1,\ldots,\#I\})=I}}{[n+1]_q \cdot \displaystyle\sum_{T \in \Lin(P)}q^{\comaj(T)}}
\end{equation}
for all $I\in \J(P)$, where
\begin{equation} \label{eqn:dist_lin_wt}
\wt(T,i) \coloneqq q^{\comaj(T,i)} \cdot q^{\#\{j \in \Des(T) \colon  j < i\}}
\end{equation}
for all $T \in \Lin(P)$ and $i=0,1,\ldots,n$.

\begin{remark} The definitions of $\vartheta(T,i)$ and $\vartheta_m(T,i)$ are quite similar. This is not a coincidence, as we will see in Corollary \ref{cor:dist_rpp_lim}.

Additionally, notice again that the exponent of $q$ in $\wt(T,i)$ is almost the comajor index of $T$ plus the number of $j \in \Des(T)$ with $j < i$, except that we always treat $i$ as a descent. This will serve as motivation for the definition of the barely set-valued comajor index for barely set-valued linear extensions used in \cref{cor:bsv_lin_expectation}.
\end{remark}

\begin{example}\label{ex:syt_probs}
In Table \ref{tab:2x2_syt} we listed the standard Young tableaux of the $\rect{2}{2}$ rectangle, and in Example \ref{ex:syt} we saw that $\displaystyle\sum_{T \in \SYT(\rect{2}{2})}q^{\comaj(T)} = q^2 + 1$. We can now compute $\mulin{q}(I)$ for all six order ideals of $\rect{2}{2}$:
\[\begin{minipage}{0.49\textwidth}
\ytableausetup{smalltableaux}
\begin{align*}
\mulin{q}\left( \varnothing \right) &= \frac{q^4(q^2+1)}{[5]_q(q^2+1)}\\
\mulin{q}\left( \parbox{0.4cm}{\centering \ydiagram{1}} \right) &= \frac{q^3(q^2+1)}{[5]_q(q^2+1)}\\
\mulin{q}\left( \parbox{0.8cm}{\centering \ydiagram{2}} \right) &= \frac{q^2}{[5]_q(q^2+1)}
\end{align*}
\end{minipage}
\begin{minipage}{0.49\textwidth}
\begin{align*}
\mulin{q}\left( \parbox{0.4cm}{\centering \ydiagram{1,1}} \right) &= \frac{q^2}{[5]_q(q^2+1)}\\
\mulin{q}\left( \parbox{0.8cm}{\centering \ydiagram{2,1}} \right) &= \frac{q(q^3+1)}{[5]_q(q^2+1)}\\
\mulin{q}\left( \parbox{0.8cm}{\centering \ydiagram{2,2}} \right) &= \frac{1\cdot (q^3+1)}{[5]_q(q^2+1)}
\end{align*}
\ytableausetup{nosmalltableaux}
\end{minipage}\]
\end{example}

The definition~\eqref{eqn:dist_lin} of $\mulin{q}$ is rather complicated, but the following proposition shows that it is indeed a distribution.

\begin{prop} \label{prop:dist_lin}
$\displaystyle \sum_{T \in \Lin(P)}\sum_{i=0}^{n} \wt(T,i) =[n+1]_q \cdot \sum_{T\in\Lin(P)} q^{\comaj(T)}$.
\end{prop}

In order to prove \cref{prop:dist_lin}, we first need a small lemma.

\begin{lemma} \label{lem:dist_lin}
Let $X\subseteq [n-1]$. Define $f_X(i)\coloneqq \#\{j \in X\colon j < i\}+\begin{cases} n-i &i \not\in X \\ 0 &i \in X\end{cases}$ for $i=0,1,\ldots,n$. Then $f_X(0)\cdots f_X(n)$ is a permutation of $0,1,\dots,n$. 
\end{lemma}
\begin{proof}
For $X=\varnothing$, we clearly have $f_\varnothing$ as the reverse permutation $w_0=n\ n-1 \dots 1\ 0$.  Let $\rho_{n,k}=(n, n-1, \dots, k)$, that is, a cycle that rotates the elements in the last $n-k+1$ positions one step right when multiplied on the right side.  We now claim that for $X=\{k_1<\cdots<k_r\}$, $f_X$ is given by the permutation $w_0\rho_{n,k_1}\dots\rho_{n,k_r}$. This is easily seen by induction. Assume it is true for $X'=\{k_1<\cdots<k_{r-1}\}$. Then clearly $f_X(i)=f_{X'}(i)$ for $i=0,\dots,k_r-1$, whereas  $f_X(k_r)=f_{X'}(k_r)-(n-k_r)$ and $f_X(i)=f_{X'}(i)+1$ for $i=k_r+1,\dots,n$. This is the same as multiplying on the right with $\rho_{n,k_r}$. An example of this induction is displayed in~\cref{fig:dist_lin}.
\end{proof}

\begin{figure}
\[\begin{array}{l | c}
 X' & f_{X'} \\ \hline 
\varnothing & 6\ 5\ 4\ 3\ 2\ 1\ 0 \\
\{2\} & 6\ 5\ 0\ 4\ 3\ 2\ 1 \\
\{2,4\} & 6\ 5\ 0\ 4\ 1\ 3\ 2 \\
\{2,4,5\} & 6\ 5\ 0\ 4\ 1\ 2\ 3 \\
\end{array} \]
\caption{An example of the proof of~\cref{lem:dist_lin}, with $n=6$ and $X=\{2,4,5\}$.} \label{fig:dist_lin}
\end{figure}

\begin{proof}[Proof of~\cref{prop:dist_lin}]
It is easy to see from~\eqref{eqn:dist_lin_wt} that
\[\wt(T,i) = q^{\comaj(T)} \cdot q^{\#\{j \in \Des(T)\colon j < i\}} \cdot q^{\begin{cases} n-i &i \not\in \Des(T) \\ 0 &i \in \Des(T)\end{cases}}\]
for any $T \in \Lin(P)$ and $i=0,1,\ldots,n$. Hence \cref{lem:dist_lin} implies that
\[ \sum_{i=0}^{n} \wt(T,i) = [n+1]_q\cdot q^{\comaj(T)}\]
for any $T \in \Lin(P)$, which proves the proposition.
\end{proof}

\begin{remark}\label{rem:lin_prob_split}
The above discussion shows that the distribution $\mulin{q}$ ``splits as a product'' in the same way that $\murpp{q}{m}$ does. That is, $\mulin{q}(I)$ can be thought of more intuitively as the probability that we select $I$ via the following process: we select $T \in \Lin(P)$ with probability proportional to $q^{\comaj(T)}$; we select $i \in \{0,1,\ldots,n\}$ independently with probability proportional to $q^{i}$; and then we select the order ideal $T^{-1}(\{1,\ldots,f_{\Des(T)}^{-1}(i)\})$, where $f_{\Des(T)}$ is the permutation of $0,1,\ldots,n$ from \cref{lem:dist_lin} corresponding to the subset $X=\Des(T)\subseteq[n-1]$.
\end{remark}

So $\mulin{q}$ indeed defines a distribution on $\J(P)$, but what is the significance of this distribution? The point is:

\begin{cor} \label{cor:dist_rpp_lim}
For $0 < q \leq 1$, we have $\displaystyle \lim_{m\to \infty} \murpp{q}{m} = \mulin{q}$ (where the limit is taken pointwise).
\end{cor}

\begin{proof}
The $q=1$ version of this statement is \cite[Proposition~2.9]{chan2017expected}, so let us assume that~$0 < q < 1$. Then,  set 
\[\displaystyle \alpha(q) \coloneqq \frac{1}{(1-q)^{n+1}\;[n+1]_q!}=\frac{1}{(1-q)(1-q^2)(1-q^3)\cdots(1-q^{n+1})}.\] 
By the definition of $\murpp{q}{m}$ in \eqref{eqn:dist_rpp} together with \cref{thm:m-weight}, we have that
\begin{equation} \label{eqn:rpp_dist_lim}
\lim_{m\to\infty} \murpp{q}{m}(I) = \frac{\displaystyle \sum_{T\in \Lin(P)} \lim_{m\to\infty} \wt_m(T,\#I)\cdot \chi(T^{-1}(\{1,\ldots,\#I\})=I)}{\displaystyle \lim_{m\to\infty} \; [m]_q \cdot \sum_{\pi\in \RPP_m(P)} q^{|\pi|}}    
\end{equation}
for all $I\in \J(P)$. For any fixed $c$ and $0<q<1$, $\lim_{m\to\infty} [m-c]_q=\sum_{i=0}^{\infty}q^i=\frac{1}{1-q}$. Thus, from \eqref{eqn:rpp_gf}, we can compute that
\begin{align} \label{eqn:rpp_dist_lim_bottom}
\displaystyle \lim_{m\to\infty} \; [m]_q \cdot \sum_{\pi\in \RPP_m(P)} q^{|\pi|} &= \frac{1}{1-q} \cdot \sum_{\pi\in \RPP(P)} q^{|\pi|} \\
\notag &= \frac{1}{1-q} \cdot \frac{\sum_{T \in \Lin(P)} q^{\comaj(T)}}{(1-q)(1-q^2)\cdots(1-q^n)} \\
\notag &= \alpha(q) \cdot [n+1]_q \cdot \sum_{T \in \Lin(P)} q^{\comaj(T)}.
\end{align}
Next, note that
\[\wt_m(T,i) = \qbinom{m+n - \#(\Des(T)\setminus\{i\})}{n+1}_q \, \wt(T,i),\]
so, letting $c\coloneqq\#(\Des(T)\setminus\{i\})$, we also have
\begin{align}  \label{eqn:rpp_dist_lim_top}
\displaystyle\lim_{m\to\infty}\wt_m(T,i) &= \lim_{m\to\infty} \qbinom{m+n - \#(\Des(T)\setminus\{i\})}{n+1}_q \, \wt(T,i) \\
\notag &= \lim_{m\to\infty}\frac{[m+n-c]_q\cdots [m-c]_q}{[n+1]_q!} \, \wt(T,i) \\
\notag &=\alpha(q) \cdot \wt(T,i).
\end{align}
Combining \eqref{eqn:rpp_dist_lim}, \eqref{eqn:rpp_dist_lim_bottom}, and~\eqref{eqn:rpp_dist_lim_top}, and recalling the definition of $\mulin{q}$ in \eqref{eqn:dist_lin}, gives
\[\lim_{m\to\infty} \murpp{q}{m}(I) = \frac{\displaystyle \sum_{T \in \Lin(P)} \alpha(q)\cdot \wt(T,i) \cdot \cdot \chi(T^{-1}(\{1,\ldots,\#I\})=I)}{\displaystyle \alpha(q) \cdot [n+1]_q \cdot \sum_{T \in \Lin(P)}q^{\comaj(T)}}=\mulin{q}(I)\]
for all $I \in \J(P)$, as claimed. 
\end{proof}

\begin{remark} \label{rem:q_gt_1_lim}
For all $I \in \J(P)$, we have $\murpp{q}{m}(I) = \murpp{q^{-1}}{m}(P\setminus I)$, where $P\setminus I$ is viewed as an order ideal of the dual poset $P^{*}$. Hence, when $q > 1$, the limit $\displaystyle \lim_{m\to \infty} \murpp{q}{m}$ can be computed using the linear extensions of $P^{*}$ instead of $P$.

Let us explain the precise formula obtained this way. There is a canonical bijection $\Lin(P)\xrightarrow{\sim}\Lin(P^{*})$ given by $T \mapsto T^*$ where $T^*(p) = n+1-T(p)$ for all $p\in P^{*}$. We have $\Des(T^*)=\{n-i\colon i \in \Des(T)\}$. So, when $q > 1$, we get for $I\in\J(P)$
\[ \displaystyle \lim_{m \to \infty} \murpp{q}{m}(I) = \frac{\sum_{T\in \Lin(P)} \wt^{*}(T,\#I) \cdot \bracket{T^{-1}(\{1,2,\ldots,\#I\})=I} }{[n+1]_{q^{-1}} \sum_{T \in \Lin(P)} q^{-\maj(T)} }\]
for all $I \in \J(P)$, where
\[ \wt^{*}(T,i) \coloneqq q^{i-n} \cdot \prod_{j\in \Des(T), j < i} q^{j-n} \cdot \prod_{j \in \Des(T), j>i}  q^{j-(n+1)}\]
for all $T \in \Lin(P)$ and $i=0,1,\ldots,n$, and the \dfn{major index} of a linear extension $T \in \Lin(P)$ is $\maj(T) \coloneqq\sum_{i \in \Des(T)}i$.
\end{remark}

\begin{remark}
Suppose that $P$ is \dfn{self-dual}, i.e., that $P\simeq P^{*}$. Then in light of \cref{rem:q_gt_1_lim}, we will have $\murpp{q}{m} = \murpp{q^{-1}}{m}$ for all $m \geq 1$. However, in general $\mulin{q} \neq \mulin{q^{-1}}$ even when $P$ is self-dual, precisely because the limit in \cref{cor:dist_rpp_lim} is only valid in the regime $0 < q \leq 1$.
\end{remark}

\subsection{The rank distribution}

We now define one last distribution on $\J(P)$ which will be useful in certain later computations (such as for the minuscule posets $\Lambda_{E_6}$ and~$\Lambda_{E_7}$); see Remark \ref{rem:rank_expectation}.

For this distribution to be defined, we need to assume that our $P$ is \dfn{graded}, i.e., that all the maximal chains of $P$ have the same length. In this case $P$ has a unique \dfn{rank function} $\rk\colon P \to \N$ satisfying:
\begin{itemize}
\item $\rk(p)=0$ if $p \in P$ is minimal;
\item $\rk(p') = \rk(p) + 1$ if $p \lessdot p' \in P$.
\end{itemize}
We define the \dfn{rank} of $P$, denoted $\rk(P)$ by abuse of notation, to be $\rk(P) \coloneqq\rk(p)$ for any maximal element $p \in P$.

We define the distribution $\murk{q}\colon \J(P)\to \R$ by
\begin{equation} \label{eqn:dist_rk} 
\murk{q}(I) = \frac{\displaystyle \sum_{i=-1}^{\rk(P)} q^{\rk(P)-i}\cdot \bracket{\rk^{-1}(\{0,1,\ldots,i\})=I} }{[\rk(P)+2]_q}
\end{equation}
for all $I\in \J(P)$. More intuitively: we select $i\in \{-1,0,\ldots,\rk(P)\}$ with probability proportional to~$q^{\rk(P)-i}$; then we select the order ideal $\rk^{-1}(\{0,1,\ldots,i\})$. 

Note that, unlike the other distributions, the support of $\murk{q}$ is in general much less than all of $\J(P)$.

\section{Toggleability statistics and \texorpdfstring{$q$}{q}-toggle-symmetry} \label{sec:tog_sym}

We continue to fix a poset $P$ and retain the notation from the previous section. In this section we identify a special class of probability distributions $\mu$ on $\J(P)$, and show that the four distributions from \cref{sec:dists} belong to this class.

For $p\in P$, define \dfn{toggleability statistics} $\Tin{p}, \Tout{p}\colon \J(P)\to\Q$ by
\begin{align*}
\Tin{p}(I) &\coloneqq\begin{cases} 1 &\textrm{if $p \in \min(P/I)$}; \\ 0 &\textrm{otherwise}, \end{cases} \\
\Tout{p}(I) &\coloneqq\begin{cases} 1 &\textrm{if $p \in \max(I)$}; \\ 0 &\textrm{otherwise}, \end{cases}
\end{align*}
for all $I\in \J(P)$. Here for a subset $S\subseteq P$ we use $\min(S)$ to denote the set of minimal elements of $S$, and similarly $\max(S)$ to denote the set of maximal elements of $S$. When $\Tin{p}(I)=1$  we say $p$ can be \dfn{toggled into} $I$; when $\Tout{p}(I)=1$  we say $p$ can be \dfn{toggled out of} $I$. This toggling terminology goes back to Striker and Williams~\cite{striker2012promotion}: \dfn{toggling} at element $p$ refers to the involutive operation on $\J(P)$ of adding $p$ to~$I$ ``if we can'' (i.e., if $p\notin I$ and $I\cup\{p\}$ remains an order ideal), removing $p$ from~$I$ if we can, and otherwise not changing $I$. Note that if $p$ can be toggled into $I$ then it cannot be toggled out of~$I$, and if $p$ can be toggled out of~$I$ then it cannot be toggled into $I$.

Now recall that we have a fixed real number parameter $q > 0$. Then define the \dfn{$q$-toggleability-statistic} $\Tq{p}{q}\colon \J(P)\to\R$ by
\[\Tq{p}{q} \coloneqq\Tin{p} - q\Tout{p}.\]
We view these statistics as random variables with respect to distributions on~$\J(P)$. Specifically, for a distribution $\mu$ on $\J(P)$, and any statistic $f\colon \J(P)\to \R$, we use the notation $\E_{\mu}(f)$ to denote the expectation $\E(X)$ of the random variable $X=f(I)$ when~$I$ is distributed according to $\mu$. That is,
\[\E_{\mu}(f) \coloneqq\sum_{I \in \J(P)} f(I) \cdot \mu(I).\]

\begin{definition}
We say that a probability distribution $\mu$ on $\J(P)$ is \dfn{$q$-toggle-symmetric} if $\E_{\mu}(\Tq{p}{q}) =0$ for all $p \in P$.
\end{definition}

In other words, to say $\mu$ is $q$-toggle-symmetric means that, for every $p\in P$, we are $q$ times as likely to be able to toggle $p$ into $I$ as to toggle $p$ out of $I$ when $I$ is distributed according to $\mu$. In the case $q=1$, this becomes the notion of toggle-symmetry that was introduced in~\cite{chan2017expected}.

We will show that all the distributions from \cref{sec:dists} are $q$-toggle-symmetric. But first, the reader may be wondering what the significance of this notion is. Without going into too much detail, the significance is that, in fortuitous situations, a statistic $f\colon \J(P)\to\R$ may have a representation
\[ f=c(q) + \sum_{p\in P} c_p(q)\Tq{p}{q}\]
for some constants $c(q), c_p(q) \in \R$ (which may depend on $q$). In this case, we know by linearity of expectations that $\E_{\mu}(f) = c(q)$ for every $q$-toggle-symmetric distribution $\mu$. The fact that the expectation of $f$ is this same constant $c(q)$ for every $q$-toggle-symmetric distribution can be used to prove combinatorial identities by double-counting. We will see this happen in \cref{sec:proof_rect}.

We now proceed to prove the $q$-toggle-symmetry of the distributions from \cref{sec:dists}.

\begin{lemma} \label{lem:rpp_tog_sym}
For any $m\geq 1$, the distribution $\murpp{q}{m}$ is $q$-toggle-symmetric.
\end{lemma}

\begin{proof}
Let $p \in P$. We need to show that $\E_{\murpp{q}{m}}(\Tq{p}{q})=0$. Equivalently, thanks to the definition~\eqref{eqn:dist_rpp} of $\murpp{q}{m}$, we need to show that
\begin{equation} \label{eqn:toggle_cancel}
 \sum_{\pi \in \RPP_m(P)} \sum_{k=0}^{m-1} q^{|\pi|+k} \cdot \Tq{p}{q}(\pi^{-1}(\{0,1,\ldots,k\})) = 0.
\end{equation}
To show~\eqref{eqn:toggle_cancel}, it suffices to define an involution $(\pi,k) \mapsto (\pi',k')$ on the set of pairs of a $\pi \in \RPP_m(P)$ and $k=0,1,\ldots,m-1$ for which
\begin{equation} \label{eqn:toggle_cancel_bij}
q^{|\pi|+k} \cdot \Tq{p}{q}(\pi^{-1}(\{0,1,\ldots,k\})) + q^{|\pi'|+k'} \cdot \Tq{p}{q}((\pi')^{-1}(\{0,1,\ldots,k'\}))=0.
\end{equation}

So now let $(\pi,k)$ be such a pair, and we will define the corresponding $(\pi',k')$. Let 
\begin{align*}
x &\coloneqq\max(\{\pi(r)\colon r < p\in P\}); \\ 
y &\coloneqq\min(\{\pi(r)\colon r > p\in P\}),
\end{align*}
with the convention that $\max(\varnothing)=0$ and $\min(\varnothing)=m$, so that $x \leq \pi(p) \leq y$. Observe that
\[\Tq{p}{q}(\pi^{-1}(\{0,1,\ldots,k\}))=\begin{cases} 1 &\textrm{if $k=x,x+1,\ldots,\pi(p)-1$}; \\ -q &\textrm{if $k=\pi(p),\pi(p)+1,\ldots,y-1$}; \\ 0&\textrm{otherwise}.\end{cases}\]
If $k < x$ or $k \geq y$, then define $(\pi',k') \coloneqq(\pi,k)$. Clearly this will satisfy~\eqref{eqn:toggle_cancel_bij} since each term in the sum will be $0$.
From now on assume that $x \leq k < y$. First suppose that $x \leq k < \pi(p)$. Then define $(\pi',k')$ by
\[ \pi'(p') \coloneqq\begin{cases} k &\textrm{if $p'=p$}; \\ \pi(p') &\textrm{otherwise} \end{cases}\]
for all $p'\in P$, and $k' \coloneqq\pi(p)-1$. We will have
\begin{align*}
q^{|\pi|+k} \cdot \Tq{p}{q}(\pi^{-1}(\{0,1,\ldots,k\})) + q^{|\pi'|+k'} \cdot \Tq{p}{q}((\pi')^{-1}(\{0,1,\ldots,k'\}))= \\
q^{|\pi|+k} \cdot 1 + q^{|\pi|+k-1} \cdot (-q) = 0,
\end{align*}
so that~\eqref{eqn:toggle_cancel_bij} is satisfied. Next suppose that $\pi(p)\leq k < y$. Then define $(\pi',k')$ by
\[ \pi'(p') \coloneqq\begin{cases} k+1 &\textrm{if $p'=p$}; \\ \pi(p) &\textrm{otherwise} \end{cases}\]
for all $p' \in P$, and $k' \coloneqq\pi(p)$. We will have
\begin{align*}
q^{|\pi|+k} \cdot \Tq{p}{q}(\pi^{-1}(\{0,1,\ldots,k\})) + q^{|\pi'|+k'} \cdot \Tq{p}{q}((\pi')^{-1}(\{0,1,\ldots,k'\}))= \\
q^{|\pi|+i} \cdot (-q) + q^{|\pi|+k+1} \cdot 1 = 0,
\end{align*}
so that~\eqref{eqn:toggle_cancel_bij} is satisfied.

It can be seen that the map $(\pi,k)\mapsto (\pi',k')$ we defined is an involution, so we are done.
\end{proof}

\begin{cor}
The distribution $\muuni{q}$ is $q$-toggle-symmetric.
\end{cor}
\begin{proof}
Since $\muuni{q}=\murpp{q}{1}$, this follows from \cref{lem:rpp_tog_sym}.
\end{proof}

\begin{cor} \label{cor:lin_tog_sym}
The distribution $\mulin{q}$ is $q$-toggle-symmetric.
\end{cor}
\begin{proof}
By \cref{cor:dist_rpp_lim}, we know that $\mulin{q}=\displaystyle \lim_{m\to\infty} \murpp{q}{m}$ for $0 < q \leq 1$. So in this regime of $q$, that  $\mulin{q}$ is $q$-toggle-symmetric follows from \cref{lem:rpp_tog_sym}. To say that $\mulin{q}$ is $q$-toggle-symmetric for all $0 < q$ is just to say that certain polynomials in $q$ are equal to zero. If they are equal to zero for $0 < q \leq 1$, then they are equal to zero for all $0 < q$.
\end{proof}

\begin{remark}
It is also possible to prove \cref{cor:lin_tog_sym} bijectively at the level of linear extensions, but it seems that the bijection needed for linear extensions is much more complicated than the one we have defined above for reverse $P$-partitions. Note that in the case $q=1$, there is a relatively straightforward bijection of linear extensions which establishes toggle-symmetry: see~\cite[Remark~2.11]{chan2017expected}.
\end{remark}

Finally, for the last distribution from \cref{sec:dists}:

\begin{lemma}
If $P$ is graded, the distribution $\murk{q}$ is $q$-toggle-symmetric.
\end{lemma}
\begin{proof}
This is straightforward from the definition~\eqref{eqn:dist_rk} of~$\murk{q}$, and the fact that for any $p\in P$,
\[ \Tq{p}{q}( \rk^{-1}(\{0,1,\ldots,i\})) = \begin{cases} 1 &\textrm{if $i=\rk(p)-1$}; \\ -q &\textrm{if $i=\rk(p)$}; \\ 0 &\textrm{otherwise},\end{cases}\]
for all $i=-1,0,\ldots,\rk(P)$.
\end{proof}

\begin{remark} \label{rem:rank_expectation}
As mentioned above, we are often interested in the expectation $\E_{\mu}(f)$ for $q$-toggle-symmetric distributions $\mu$ and statistics $f$ which happen to have the form
\[ f=c(q) + \sum_{p\in P} c_p(q)\Tq{p}{q}.\]
We see that $\E_{\mu}(f)=c(q)$. Suppose we were not given the constant~$c(q)$ ahead of time. In order to compute $c(q)$, it is most convenient (in the case when $P$ is graded) to take $\mu=\murk{q}$. This is because $\murk{q}$ has sparse support. Explicitly, we get
\[ c(q) = \E_{\murk{q}}(f)=\frac{1}{[\rk(P)+2]_q}\sum_{i=-1}^{\rk(P)} q^{\rk(P)-i} \cdot f(\rk^{-1}(\{0,1,\ldots,i\}) ).\]
\end{remark}

\section{Barely set-valued fillings and proofs of the main results} \label{sec:proof_rect}

In this section we connect the discussion of distributions and expectation to the focus of \cref{sec:intro}: set-valued fillings, and more specifically, \emph{barely} set-valued fillings. We continue to work with an arbitrary poset $P$ on $n$ elements as above; we will only specialize to the rectangle at the last moment.

We recall that set-valued tableaux were introduced by Buch~\cite{buch2002littlewood} in his study of the $K$-theory of Grassmannians. Set-valued $P$-partitions were introduced by Lam and Pylyavksyy~\cite{lam2007combinatorial} in their development of a $K$-theoretic version of the theory of quasisymmetric functions. 

A \dfn{set-valued filling} $\tau\colon P \to 2^{\N}$ of $P$ is a map which assigns to each element of $P$ a finite, non-empty set of natural numbers: $\varnothing \neq \tau(p) \subset_{\mathrm{fin}} \N$ for all $p \in P$. We say that the set-valued filling $\tau$ is a \dfn{set-valued reverse $P$-partition} if $\max(\tau(p)) \leq \min(\tau(q))$ for all $p < q \in P$. If all entries are singletons, then such a $\tau$ is really just a reverse $P$-partition. We say that the set-valued reverse $P$-partition $\tau$ is \dfn{barely set-valued} if there is a (necessarily unique) $\pstar(\tau)\in P$ for which:
\begin{itemize}
\item $\#\tau(p) = 1$ for all $p\neq \pstar(\tau)\in P$;
\item $\#\tau(\pstar(\tau)) = 2$.
\end{itemize}
We use $\RPPo(P)$ to denote the set of barely set-valued reverse $P$-partitions and use $\RPPo_m(P)$ to denote the set of $\tau \in \RPPo(P)$ with $\tau(p)\subseteq \{0,1,\ldots,m\}$ for all $p \in P$. For $\tau \in \RPPo(P)$, we use the notation $\istar(\tau) \coloneqq\max(\tau(\pstar(\tau)))$.

Following Reiner--Tenner--Yong~\cite{reiner2018poset}, we observe that $\RPPo_m(P)$ has a simple description in terms of $\RPP_m(P)$:

\begin{lemma} \label{lem:bsv_rpp_bij}
For any $m \geq 1$, we have a bijection
\begin{align*}
\RPPo_m(P) &\xrightarrow{\sim} \left\{ (\pi,i,p)\colon \parbox{2.25in}{\begin{center} $\pi \in \RPP_m(P)$,\\ $i=0,1,\ldots,m-1,$ \\ $p \in \max(\pi^{-1}(\{0,1,\ldots,i\}))$ \end{center}} \right\} \\
\tau &\mapsto (\pi,\, \istar(\tau)-1, \, \pstar(\tau))
\end{align*}
where $\pi$ is defined by
\[ \pi(p) \coloneqq\min(\tau(p))\]
for all $p\in P$.
\end{lemma}

\ytableausetup{nosmalltableaux}
\begin{example}
Consider the following barely set-valued reverse plane partition $\tau \in \RPP_{6}^{+1}(\rect{3}{3})$: 
\[ \begin{ytableau}
    1 & 1 & 2\\
    3 &*(yellow) 4, 5 & 5\\
    3 & 5 & 6
    \end{ytableau}\]
This corresponds to the triple $(\pi, 4, (2,2))$, where $\pi$ is the (ordinary) reverse plane partition shown below:
\[\begin{ytableau}
1 & 1 & 2\\
3 & 4 & 5\\
3 & 5 & 6
\end{ytableau}\]
\end{example}

\begin{proof} [Proof of \cref{lem:bsv_rpp_bij}]
The inverse map is defined as follows. Given a triple $(\pi,i,p)$, we define a $\tau \in \RPP_m^{+1}(P)$ by
\[\tau(p') = \begin{cases}\{\pi(p')\}, & p' \neq p\\ \{\pi(p'),i+1\}, & p'=p\end{cases}.\]
Since by definition $p \in \max(\pi^{-1}(\{0,\dots,i\}))$, we know that $\tau(p') \geq i+1$ for all $p'> p$. Thus $\tau$ is indeed a reverse $P$-partition, and since there is exactly one $p$ for which $\#\tau(p) = 2$, $\tau$ is barely set-valued.
\end{proof}

Thanks to \cref{lem:bsv_rpp_bij}, we can count (and even $q$-count) $\RPPo_m(P)$ using expectations of a certain random variable, as briefly suggested at the end of \cref{sec:intro}. Specifically, define the \dfn{down-degree} statistic $\ddeg\colon \J(P)\to\Q$ by 
\[\ddeg(I) \coloneqq\#\max(I)\] 
for all $I\in \J(P)$. Note that equivalently we have
\[ \ddeg = \sum_{p\in P}\Tout{p}.\]
This statistic is called down-degree because it is the down-degree of (i.e., number of elements covered by) the order ideal $I$ in the distributive lattice $\J(P)$.
The statistic with respect to which we want to $q$-count $\RPPo_m(P)$ is very easy to define: for $\tau \in \RPPo_m(P)$, we set $|\tau| \coloneqq\sum_{p\in P} \sum_{i \in \tau(p)} i$.

\begin{cor} \label{cor:bsv_rpp_expectation}
For any $q > 0$ and $m\geq 1$, we have
\[\E_{\murpp{q}{m}}(\ddeg) = \frac{\displaystyle \sum_{\tau\in \RPPo_m(P)}q^{|\tau|-1} }{\displaystyle [m]_q \cdot \sum_{\pi\in \RPP_m(P)}q^{|\pi|}}.\]
\end{cor}
\begin{proof}
Under the bijection in \cref{lem:bsv_rpp_bij}, we have
\[\sum_{\tau \in \RPP_m^{+1}(P)}q^{|\tau| - 1} = \sum_{(\pi,i,p)} q^{|\pi|+i}.\]
On the right-hand side of this equality, for any fixed values of $\pi$ and $i$ we obtain the same order ideal $I \in \J(P)$ as $\pi^{-1}(\{0,\dots,i\})$ for all values of $p$; hence, this ideal occurs with multiplicity equal to $\#\max(I) = \ddeg(I)$. This means that we can rewrite the right-hand side as
\[\sum_{I \in \J(P)}\left(\sum_{\pi \in \RPP_m(P)}\sum_{i=0}^{m-1}q^{|\pi|+i}\cdot\bracket{\pi^{-1}(\{0,\dots,i\}) = I}\right)\cdot \ddeg(I).\]
By definition~\eqref{eqn:dist_rpp} we thus get
\begin{align*}
\E_{\murpp{q}{m}}(\ddeg)= \sum_{I \in \J(P)}\murpp{q}{m}(I)\cdot\ddeg(I)=\frac{\displaystyle \sum_{\tau\in \RPPo_m(P)}q^{|\tau|-1} }{\displaystyle [m]_q \cdot \sum_{\pi\in \RPP_m(P)}q^{|\pi|}}.
\end{align*}
\end{proof}

\begin{example}
In Example \ref{ex:pp} we computed the two sums on the top and bottom of the right-hand side of the equation in~\cref{cor:bsv_rpp_expectation} for $m=1$ and $P=\rect{2}{2}$. In Example \ref{ex:rpp_probs} we computed~$\murpp{q}{1}(I)$ for all six order ideals of $\rect{2}{2}$. 
Combining that and the definition of $\ddeg$ gives us 
\[\E_{\murpp{q}{1}}(\ddeg) = \frac{q^3 + 2q^2 + 2q + 1}{q^4+q^3+2q^2 + q + 1}.\]
Since $[m]_q=1$ in this case, we see that \cref{cor:bsv_rpp_expectation} indeed holds in this case.
\end{example}

We can of course also consider set-valued analogs of linear extensions instead of reverse $P$-partitions. We define a \dfn{set-valued linear extension} of $P$ to be a set-valued filling $S\colon P \to 2^{\N}$ for which:
\begin{itemize}
\item $\bigcup_{p\in P} S(p) = [n+k]$ for some $k \geq 0$;
\item $S(p) \cap S(q) = \varnothing$ for $p\neq q \in P$;
\item $\max(S(p)) < \min(S(q))$ for $p < q \in P$.
\end{itemize}
In other words, the first two items say that $\{S(p)\colon p \in P\}$ is a partition of the set $[n+k]$, and the third item says that the filling is order-preserving. In the case $k=0$ such a set-valued linear extension is really just a linear extension. In the case $k=1$ we say that such a $S$ is \dfn{barely set-valued}. We use $\Lino(P)$ to denote the set of barely set-valued linear extensions of $P$. For $S \in \Lino(P)$ there is a unique $\pstar(S)\in P$ with $\#S(\pstar(S))=2$; we also use the notation $\istar(S) \coloneqq\max(S(\pstar(S)))$.

Again, there is a simple description of $\Lino(P)$ in terms of $\Lin(P)$:

\begin{lemma} \label{lem:bsv_lin_bij}
We have a bijection
\begin{align*}
\Lino(P) &\xrightarrow{\sim} \left\{ (T,i,p)\colon \parbox{2.25in}{\begin{center} 
$T \in \Lin(P), i=0,1,\ldots,n,$ \\ $p \in \max(T^{-1}(\{1,\ldots,i\}))$ \end{center}} \right\} \\
S &\mapsto (T, \,  \istar(S)-1, \,  \pstar(S))
\end{align*}
where $T$ is defined by
\[ T(p) \coloneqq\begin{cases} \min(S(p)) &\textrm{if $\min(S(p)) < \istar(S)$}; \\ \min(S(p))-1 &\textrm{otherwise},\end{cases}\]
for all $p\in P$.
\end{lemma}

\begin{example}
Consider the barely set-valued linear extension $S \in \Lino(\rect{3}{3})$ represented as a barely set-valued standard Young tableau below:
\[\begin{ytableau}
1 & 2 & 5\\
3 & *(yellow)4, 6 & 8\\
7 & 9 & 10
\end{ytableau}\]
Under the bijection above, $S$ corresponds to the triple $(T, 5, (2,2))$, where $T \in \Lin(\rect{3}{3})$ is the linear extension represented as the following standard Young tableau:
\[\begin{ytableau}
1 & 2 & 5\\
3 & 4 & 7\\
6 & 8 & 9
\end{ytableau}\]
\end{example}
\begin{proof} [Proof of \cref{lem:bsv_lin_bij}]
The inverse map is defined as follows. Given a tuple $(T,i,p)$ we define a barely set-valued function $S$ on $P$ by
\[S(p') = \begin{cases} \{T(p')\}, & p'\neq p, \; T(p')\leq i\\ \{T(p') + 1\}, & p'\neq p, \; T(p') >i\\ \{T(p), i+1\}, & p'=p \end{cases}.\]
We need to check that $S$ is a barely set-valued linear extension. Since $T$ is a linear extension we see that $\min S(p') < \min S(p'')$ whenever $p' < p''$. Moreover, we know that if $p' > p$ then $S(p') > i+1$, so we see also that $\max S(p) < \min S(p')$ whenever~$p' > p$. 
\end{proof}

As before, we can use \cref{lem:bsv_lin_bij} to $q$-count $\Lino(P)$, although now the statistic we need is more involved.  Let $S \in \Lino(P)$ and $\omega$ a natural labeling of $P$. A \dfn{barely set-valued (BSV) descent} of $S$ is a number $i=1,2,\ldots,n+1$ with $\omega(p') < \omega(p)$, where $p',p\in P$ are such that $i+1 \in S(p')$ and $i \in S(p)$, \emph{except that}:
\begin{itemize}
\item $\istar(S)-1$ is \emph{never} a BSV descent;
\item $\istar(S)$ is \emph{always} a BSV descent.
\end{itemize}
We denote the set of BSV descents of $S$ by $\Des^{+1}(S)$. We define the \dfn{barely set-valued comajor index} of $S$ to be $\comaj^{+1}(S) \coloneqq\sum_{i \in \Des^{+1}(S)} (n+1-i)$.

\begin{cor} \label{cor:bsv_lin_expectation}
For any $q > 0$, we have
\[\E_{\mulin{q}}(\ddeg) = \frac{\displaystyle \sum_{S\in \Lino(P)}q^{\comaj^{+1}(S)} }{\displaystyle [n+1]_q \cdot \sum_{T\in \Lin(P)}q^{\comaj(T)}}.\]
\end{cor}
\begin{proof}
This follows from the definition~\eqref{eqn:dist_lin} of the distribution $\mulin{q}$, and the bijection in \cref{lem:bsv_lin_bij} in the same way that \cref{cor:bsv_rpp_expectation} follows from~\eqref{eqn:dist_rpp} and the bijection in \cref{lem:bsv_rpp_bij}.
\end{proof}

\begin{example}
In Example \ref{ex:syt_probs} we computed $\mulin{q}(I)$ for the six order ideals of $\rect{2}{2}$. 
Using these computations and the definition of $\ddeg$, we get
\begin{align*}\E_{\mulin{q}}(\ddeg) &= \frac{q^3(q^2+1) +  2q^2 + 2q(q^3+1) + q^3 + 1}{[5]_q(q^2+1)}\\
&= \frac{q^5 + 2q^4 + 2q^3 +  2q^2 + 2q + 1}{[5]_q(q^2+1)}.\end{align*}
In Example \ref{ex:syt}, we saw that $\sum_{S \in \Lino(\rect{2}{2})}q^{\comaj^{+1}(S)}=q^5 + 2q^4 + 2q^3 + 2q^2 + 2q + 1$, as expected by \cref{cor:bsv_lin_expectation}.
\end{example}

Even with \cref{cor:bsv_rpp_expectation,cor:bsv_lin_expectation}, it is not obvious how we arrive at product formulas for the $q$-enumerations of $\RPPo_m(P)$ and $\Lino(P)$. Indeed, for general posets $P$ we have no hope of deriving such formulas. However, certain special posets have the miraculous property that these expectations have a simple form that can be computed using the $q$-toggle-symmetry approach.

One such special case is the rectangle poset $\rect{a}{b}$. We appeal to the following key result, which was proved in~\cite{defant2021homomesy}.

\begin{thm}[{Defant et al.~\cite{defant2021homomesy}}] \label{thm:ddeg_expression}
For $P=\rect{a}{b}$, there are constants $c_p(q) \in \Q(q)$ for which
\[ \ddeg = \frac{[a]_q[b]_q}{[a+b]_q} + \sum_{p\in P} c_p(q) \Tq{p}{q}\]
as an equality of functions $\J(P)\to \Q$. (Here we view $q$ as a formal indeterminate.)
\end{thm}

We note that the case $q=1$ of \cref{thm:ddeg_expression} was established in~\cite{chan2017expected}, and that the proof of \cref{thm:ddeg_expression} in~\cite{defant2021homomesy} is an extension of the techniques of~\cite{chan2017expected}.

\Cref{thm:ddeg_expression} now easily allows us to give product formulas for the $q$-enumerations of $\RPPo_m(\rect{a}{b})$ and $\Lino(\rect{a}{b})$.

\begin{cor} \label{cor:bsv_rpp}
For $P=\rect{a}{b}$, and any $m \geq 1$,
\begin{align*}
\sum_{\tau\in \RPPo_m(P)}q^{|\tau|-1} &= \frac{[a]_q[b]_q}{[a+b]_q} \cdot [m]_q \cdot \sum_{\pi \in \RPP_m(P)}q^{|\pi|}\\
&= \frac{[a]_q[b]_q}{[a+b]_q} \cdot [m]_q \cdot \prod_{i=1}^{a}\prod_{j=1}^{b} \frac{[i+j+m-1]_q}{[i+j-1]_q}.
\end{align*}
\end{cor}

\begin{cor} \label{cor:bsv_lin}
For $P=\rect{a}{b}$,
\begin{align*}
\sum_{S\in \Lino(P)}q^{\comaj^{+1}(S)} &= \frac{[a]_q[b]_q}{[a+b]_q} \cdot [ab+1]_q \cdot \sum_{T \in \Lin(P)}q^{\comaj(T)}\\
&= \frac{[a]_q[b]_q}{[a+b]_q} \cdot [ab+1]_q! \cdot \prod_{i=0}^{a-1} \frac{[i]_q!}{[b+i]_q!}.
\end{align*} 
\end{cor}

\begin{proof}[Proofs of \cref{cor:bsv_rpp,cor:bsv_lin}]
These corollaries follow from double-counting arguments. Since we know that $\murpp{q}{m}$ and $\mulin{q}$ are $q$-toggle-symmetric (see \cref{lem:rpp_tog_sym,cor:lin_tog_sym}), we have that $\E_{\mulin{q}}(\ddeg) = \E_{\murpp{q}{m}}(\ddeg) = \frac{[a]_q[b]_q}{[a+b]_q}$ by \cref{thm:ddeg_expression}. We substitute this expectation into the left-hand sides of \cref{cor:bsv_rpp_expectation} and \cref{cor:bsv_lin_expectation}, and then clear denominators. To get the final equalities, we use the known product formulas for $\sum_{\pi \in \RPP_m(a\times b)}q^{|\pi|}$ and $\sum_{T \in \Lin(a\times b)}q^{\comaj(T)}$, i.e., MacMahon's formula~\eqref{eqn:macmahon} and the q-hook length formula~\eqref{eqn:qhlf}.
\end{proof}

\Cref{cor:bsv_lin} directly translates to \cref{thm:syt_intro} from \cref{sec:intro}. Similarly, \cref{cor:bsv_rpp} translates to \cref{thm:pp_intro} from \cref{sec:intro}: for $P=\rect{a}{b}$, a reverse $P$-partition is almost the same thing as an $\rect{a}{b}$ plane partition, except that the former increases along rows and columns while the latter decreases along rows and columns. However, note that $P$ is self-dual, and if we take the dual of a reverse $P$-partition (i.e., rotate it by $180^\circ$), we get a plane partition.

Recall the $q,t$-number notation $[k]_{q,t} \coloneqq  \frac{(t^k-q^k)}{(t-q)} = t^{k-1}+qt^{k-2}+\cdots+q^{k-2}t+q^{k-1}$. It was mentioned in \cref{rem:rows} that we can prove the following refinement of \cref{cor:bsv_lin}.

\begin{cor} \label{cor:bsv_lin_refined}
For $P=\rect{a}{b}$,
\[\sum_{S\in \Lino(P)}q^{\comaj^{+1}(S)} t^{\rstar(S)-1} =  \frac{[a]_{q,t}[b]_q}{[a+b]_q} \cdot [ab+1]_q \cdot \sum_{T\in \Lin(P)}q^{\comaj(T)},\]
where $\rstar(S)$ denotes the row that $\pstar(S)$ occurs in.
\end{cor}

Indeed, the proof of \cref{cor:bsv_lin_refined} is exactly the same as the proof of \cref{cor:bsv_lin}, but uses the fact that for each $i=1,2,\ldots,a$ we have
\begin{equation} \label{eqn:ddeg_expression_refined}
 \sum_{j=1}^{b} \Tout{(i,j)} = \frac{q^{a-i} \cdot [b]_q}{[a+b]_q} + \sum_{p\in P} c_p(q) \Tq{p}{q}
\end{equation}
for certain constants $c_p(q) \in \Q(q)$. Observe how \eqref{eqn:ddeg_expression_refined} refines \cref{thm:ddeg_expression}. This refinement~\eqref{eqn:ddeg_expression_refined} of \cref{thm:ddeg_expression} was also proved in~\cite{defant2021homomesy}. 

In the same way we can also prove the following.

\begin{cor} \label{cor:bsv_rpp_refined}
For $P=\rect{a}{b}$, and any $m\geq 1$,
\[\sum_{\tau\in \RPPo_m(P)}q^{|\tau|-1} t^{\rstar(\tau)-1} =  \frac{[a]_{q,t}[b]_q}{[a+b]_q} \cdot [m]_q \cdot \sum_{\pi\in \RPP_m(P)}q^{|\pi|},\]
where $\rstar(\tau)$ denotes the row that $\pstar(\tau)$ occurs in.
\end{cor}

Via the transposition symmetry of the rectangle, there are of course versions of \cref{cor:bsv_lin_refined,cor:bsv_rpp_refined} which keep track of the column in which the double entry occurs instead of the row.

\section{Shifted staircase and other minuscule posets} \label{sec:sstair}

In this section we extend our results on the rectangle to shifted shapes, and beyond.

Let $\lambda$ be a \dfn{strict} partition of $n$, i.e., a partition whose parts strictly decrease. The \dfn{shifted Young diagram} of $\lambda$ is like the usual Young diagram, except that each row is indented by one box compared to the row above it. 
We continue to use the same matrix coordinates for boxes of shifted Young diagrams, so that $(i,i)$ is always the leftmost box in the $i$th row. A \dfn{(shifted) standard Young tableau} of shifted shape $\lambda$ is a filling of the shifted Young diagram of $\lambda$ with the numbers $1,\ldots,n$, each appearing once, so that entries are increasing along rows and down columns. The shifted Young diagram corresponding to $\lambda=(4,3,1)$ and an example of a shifted SYT of shape $\lambda$ appear below.
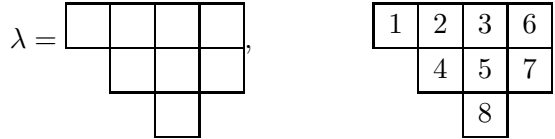
\[\lambda = \ydiagram{4,1+3,2+1},\qquad \qquad\begin{ytableau} 1 & 2 & 3 & 6 \\ \none & 4 & 5 &7\\ \none & \none & 8 \\ \none  \end{ytableau}\]

To any strict partition $\lambda$ we associate a poset, denoted $\sshape{\lambda}$, whose elements are the boxes of the shifted Young diagram of $\lambda$ and with $u\leq u'$ if $u$ is weakly northwest of $u'$. We often refer to posets of the form $\sshape{\lambda}$ as \dfn{shifted shapes}. Clearly the SYTs of shifted shape $\lambda$ are then the linear extensions of $\sshape{\lambda}$. We always assume that $\sshape{\lambda}$ is given the natural labeling $\omega$ where the 1st row of~$\lambda$ gets the values $1,2,\ldots,\lambda_1$, the 2nd row gets the values $\lambda_1+1,\ldots,\lambda_2$, and so on. Hence, a \dfn{(natural) descent} of a shifted SYT is again an entry $i$ such that $i+1$ appears in a higher row than $i$.

There is also a hook length formula for shifted shapes: for any strict partition~$\lambda$ of $n$ we have
\begin{equation} \label{eqn:hlf_shifted}
\#\Lin(\sshape{\lambda}) = n! \cdot \prod_{u \in \sshape{\lambda}} \frac{1}{h^*(u)},
\end{equation}
where $h^*(u)$ is the \dfn{shifted hook length} of the box $u$. The formula~\eqref{eqn:hlf_shifted} is usually attributed to Thrall~\cite{thrall1952combinatorial}, although see~\cite{sagan1990ubiquitous} for a more precise account of the history, as well as, e.g., the definition of $h^*(u)$.

The $q$-version of the shifted hook length formula is as follows. Gansner~\cite{gansner1978matrix} proved that
\begin{equation} \label{eqn:rpp_shifted}
\sum_{\pi \in \RPP(\sshape{\lambda})} q^{|\pi|} = \prod_{u \in \sshape{\lambda}} \frac{1}{(1-q^{h^*(u)})}.
\end{equation}
Combining~\eqref{eqn:rpp_shifted} and~\eqref{eqn:rpp_gf} gives
\begin{equation} \label{eqn:qhlf_shifted}
\sum_{T\in\Lin(\sshape{\lambda})}q^{\comaj(T)} = [n]_q! \cdot \prod_{u \in \sshape{\lambda}} \frac{1}{[h^*(u)]_q}.
\end{equation}

The formula~\eqref{eqn:qhlf_shifted} applies to all shifted shapes. However, just as in \cref{sec:proof_rect} we were ultimately only able to give product formulas for barely set-valued fillings of rectangular shapes, here we will only be able to give formulas for barely set-valued fillings of a very special shifted shape. This special shifted shape is the \dfn{shifted staircase} $\sstair{k} \coloneqq (k,k-1,\ldots,1)$. By abuse of notation we also use $\sstair{k}$ to mean the associated poset $\sstair{k} \coloneqq \sshape{\sstair{k}}$.

Recall that one thing which is special about the rectangle, compared to other shapes, is that there is not just a product formula for $\sum_{\pi \in \RPP(\rect{a}{b})} q^{|\pi|}$, there is in fact a product formula for $\sum_{\pi \in \RPP_m(\rect{a}{b})} q^{|\pi|}$ for any $m\geq 0$: namely, MacMahon's formula~\eqref{eqn:macmahon}. The same is true of the shifted staircase:
\begin{equation} \label{eqn:rpp_shifted_staircase}
\sum_{\pi \in \RPP_m(\sstair{k})} q^{|\pi|} = \prod_{1\leq i \leq j \leq k} \frac{[i+j+m-1]_q}{[i+j-1]_q}
\end{equation}
The formula~\eqref{eqn:rpp_shifted_staircase} was conjectured by Bender and Knuth~\cite{bender1972enumeration}, and proved independently by Gordon~\cite{gordon1983bender}, Andrews~\cite{andrews1977plane}, and Macdonald~\cite{macdonald1979symmetric}.

Our techniques for counting/$q$-counting barely set-valued fillings also apply to the shifted staircase. This is because Defant et al.~also proved:

\begin{thm}[{Defant et al.~\cite{defant2021homomesy}}] \label{thm:shifted_ddeg_expression}
For $P=\sstair{k}$, there are constants $c_p(q) \in \Q(q)$ for which
\[ \ddeg = \frac{\qbinom{k+1}{2}_q}{[2k]_q} + \sum_{p\in P} c_p(q) \Tq{p}{q}.\]
\end{thm}

The case $q=1$ of \cref{thm:shifted_ddeg_expression} was proved in \cite{hopkins2017CDE}.

\begin{cor} \label{cor:shifted_bsv_rpp}
For $P=\sstair{k}$, and any $m \geq 1$,
\[\sum_{\tau\in \RPPo_m(P)}q^{|\tau|-1} =   \frac{\qbinom{k+1}{2}_q}{[2k]_q} \cdot [m]_q \cdot \sum_{\pi\in \RPP_m(P)}q^{|\pi|}.\] 
\end{cor}

\begin{cor} \label{cor:shifted_bsv_lin}
For $P=\sstair{k}$,
\[\sum_{S\in \Lino(P)}q^{\comaj^{+1}(S)} = \frac{\qbinom{k+1}{2}_q}{[2k]_q} \cdot [k(k+1)/2+1]_q \cdot \sum_{T\in \Lin(P)}q^{\comaj(T)}.\] 
\end{cor}

\begin{proof}[Proofs of \cref{cor:shifted_bsv_rpp,cor:shifted_bsv_lin}]
As before, these follow from combining \cref{thm:shifted_ddeg_expression} with \cref{cor:bsv_rpp_expectation,cor:bsv_lin_expectation}.
\end{proof}

Notice that, in light of \eqref{eqn:rpp_shifted_staircase} and \eqref{eqn:qhlf_shifted}, \cref{cor:shifted_bsv_rpp,cor:shifted_bsv_lin} indeed yield product formulas for the generating functions of barely set-valued plane partitions and tableaux of shifted staircase shape.

\begin{table}
\renewcommand{\arraystretch}{1.5}
\begin{tabular} {c | c | c}
$T$ & \parbox{0.75in}{\begin{center}\begin{ytableau} 1 & 2 & 3 \\ \none & 4 & 5 \\ \none & \none & 6 \end{ytableau}\end{center}} & \parbox{0.75in}{\begin{center}\begin{ytableau} 1 & 2 & 4 \\ \none & 3 & 5 \\ \none & \none & 6 \end{ytableau}\end{center}} \\ \hline
$\Des(T)$ & $\varnothing$ & $\{3\}$ \\ \hline
$\comaj(T)$ & $0$ & $3$
\end{tabular} 
\medskip
\caption{The tableaux in $\Lin(\sstair{3})$ with (natural) descents and comajor indices.} \label{tab:d3_syt}
\end{table}

\begin{table}
\renewcommand{\arraystretch}{1.5}
\begin{tabular} {c | c | c | c | c | c}
$S$ & \parbox{0.75in}{\begin{center}\begin{ytableau} 1,2 & 3 & 4 \\ \none & 5 & 6 \\ \none & \none & 7 \end{ytableau}\end{center}} & \parbox{0.75in}{\begin{center}\begin{ytableau} 1,2 & 3 & 5 \\ \none & 4 & 6 \\ \none & \none & 7 \end{ytableau}\end{center}} & \parbox{0.75in}{\begin{center}\begin{ytableau} 1 & 2,3 & 4 \\ \none & 5 & 6 \\ \none & \none & 7 \end{ytableau}\end{center}} & \parbox{0.75in}{\begin{center}\begin{ytableau} 1 & 2,3 & 5 \\ \none & 4 & 6 \\ \none & \none & 7 \end{ytableau}\end{center}} & \parbox{0.75in}{\begin{center}\begin{ytableau} 1 & 2 & 3,4 \\ \none & 5 & 6 \\ \none & \none & 7 \end{ytableau}\end{center}} \\ \hline
$\Des(S)$ & $\{2\}$ & $\{2,4\}$ & $\{3\}$ & $\{3,4\}$ & $\{4\}$ \\ \hline
$\comaj^{+1}(S)$ & $5$ & $8$ & $4$ & $7$ & $3$
\end{tabular} 
\medskip \hrule \medskip
\begin{tabular} {c | c | c | c | c | c}
$S$ & \parbox{0.75in}{\begin{center}\begin{ytableau} 1 & 2 & 3,5 \\ \none & 4 & 6 \\ \none & \none & 7 \end{ytableau}\end{center}} & \parbox{0.75in}{\begin{center}\begin{ytableau} 1 & 2 & 4,5 \\ \none & 3 & 6 \\ \none & \none & 7 \end{ytableau}\end{center}} & \parbox{0.75in}{\begin{center}\begin{ytableau} 1 & 2 & 5 \\ \none & 3,4 & 6 \\ \none & \none & 7 \end{ytableau}\end{center}} & \parbox{0.75in}{\begin{center}\begin{ytableau} 1 & 2 & 4 \\ \none & 3,5 & 6 \\ \none & \none & 7 \end{ytableau}\end{center}} & \parbox{0.75in}{\begin{center}\begin{ytableau} 1 & 2 & 3 \\ \none & 4,5 & 6 \\ \none & \none & 7 \end{ytableau}\end{center}} \\ \hline
$\Des(S)$ & $\{5\}$ & $\{3,5\}$ & $\{4\}$ & $\{3,5\}$ & $\{5\}$ \\ \hline
$\comaj^{+1}(S)$ & $2$ & $6$ & $3$ & $6$ & $2$
\end{tabular} 
\medskip \hrule \medskip
\begin{tabular} {c | c | c | c | c }
$S$ & \parbox{0.75in}{\begin{center}\begin{ytableau} 1 & 2 & 3 \\ \none & 4 & 5,6 \\ \none & \none & 7 \end{ytableau}\end{center}} & \parbox{0.75in}{\begin{center}\begin{ytableau} 1 & 2 & 4 \\ \none & 3 & 5,6 \\ \none & \none & 7  \end{ytableau}\end{center}} & \parbox{0.75in}{\begin{center}\begin{ytableau} 1 & 2 & 3 \\ \none & 4 & 5 \\ \none & \none & 6,7  \end{ytableau}\end{center}} & \parbox{0.75in}{\begin{center}\begin{ytableau} 1 & 2 & 4 \\ \none & 3 & 5 \\ \none & \none & 6,7  \end{ytableau}\end{center}}  \\ \hline
$\Des(S)$ & $\{6\}$ & $\{3,6\}$ & $\{7\}$ & $\{3,7\}$  \\ \hline
$\comaj^{+1}(S)$ & $1$ & $5$ & $0$ & $4$ 
\end{tabular} 
\medskip
\caption{The tableaux in $\Lino(\sstair{3})$ with barely set-valued descents and comajor indices.} \label{tab:d3_bsv_syt}
\end{table}

\begin{example}
Consider the case $k=3$. There are two elements of $\Lin(\sstair{3})$, which we show together with their natural descent sets and comajor indices in~\cref{tab:d3_syt}. The generating function of their comajor indices is
\[ \sum_{T\in \Lin(\sstair{3})} q^{\comaj(T)} = q^3 + 1 = \frac{[6]_q!}{[5]_q[4]_q[3]_q[3]_q[2]_q[1]_q},\]
in agreement with~\eqref{eqn:qhlf_shifted}. Meanwhile, there are $14$ elements of $\Lino(\sstair{3})$, which together with their descent sets and comajor indices are shown in~\cref{tab:d3_bsv_syt}. The generating function of their comajor indices is
\[ \sum_{S\in \Lino(\sstair{3})} q^{\comaj^{+1}(S)} =  q^8 + q^7 + 2 q^6 + 2 q^5 + 2q^4 + 2 q^3 + 2 q^2 + q + 1 =  \frac{\qbinom{4}{2}_q}{[6]_q} \cdot  [7]_q \cdot (q^3 + 1),\]
in agreement with \cref{cor:shifted_bsv_lin}.
\end{example}

\begin{remark}
In~\cite{defant2021homomesy}, Defant et al.~proved that, for $P=\sstair{k}$, there are constants $c_p(q)\in\Q(q)$ for which
\[ \sum_{i=1}^{k} \Tout{(i,i)} = \frac{1}{1+q} + \sum_{p\in P} c_p(q) \Tq{p}{q}.\]
In the same way that we were able to obtain the refinement~\cref{cor:bsv_lin_refined} of \cref{cor:bsv_lin}, we can thus obtain the following refinement of~\cref{cor:shifted_bsv_lin}:
\[\sum_{S\in \Lino(\sstair{k})} q^{\comaj^{+1}(S)}t^{\dstar(S)} = \frac{q\qbinom{k}{2}_q + t[k]_{q^2}}{[2k]_q} \cdot [k(k+1)/2 + 1]_q \cdot  \sum_{T\in \Lin(\sstair{k})} q^{\comaj(T)}.\]
Here for $S\in \Lino(\sstair{k})$, $\dstar(S)$ is one if the box $\pstar(S)$ with a double entry is on the main diagonal (i.e., equals $(i,i)$ for some $i$), and is zero otherwise.
\end{remark}

\begin{figure}
\begin{center}
\begin{tikzpicture}[scale=0.4]
	\node[shape=circle,fill=black,inner sep=1.5] (1) at (0,0) {};
	\node[shape=circle,fill=black,inner sep=1.5] (2) at (-1,1) {};
	\node[shape=circle,fill=black,inner sep=1.5] (3) at (-2,2) {};
	\node[shape=circle,fill=black,inner sep=1.5] (4) at (-3,3) {};
	\node[shape=circle,fill=black,inner sep=1.5] (5) at (-1,3) {};
	\node[shape=circle,fill=black,inner sep=1.5] (6) at (-2,4) {};
	\node[shape=circle,fill=black,inner sep=1.5] (7) at (-3,5) {};
	\node[shape=circle,fill=black,inner sep=1.5] (8) at (-4,6) {};
	\draw [thick] (1) -- (2) -- (3) -- (4) -- (6) -- (7) -- (8);
	\draw [thick] (3) -- (5) -- (6);
	\draw [decorate,decoration={brace,amplitude=10pt},yshift=8.5pt] (-4,6) -- (-0.9,2.9) node [black,midway,yshift=0.4cm,xshift=0.4cm]  {\rotatebox{-45}{\footnotesize $k$}};
	\rotatebox{-45}{\draw [decorate,decoration={brace,mirror,amplitude=10pt},yshift=-8.5pt] (-4.2,0) -- (0,0) node [black,midway,yshift=-0.6cm]  {\footnotesize $k$};}
	\node at (-2,-4) {$D(k)$};
\end{tikzpicture} \qquad \vline \qquad
\begin{tikzpicture}[scale=0.4]
	\node[shape=circle,fill=black,inner sep=1.5] (-4) at (3,-3) {};
	\node[shape=circle,fill=black,inner sep=1.5] (-3) at (2,-2) {};
	\node[shape=circle,fill=black,inner sep=1.5] (-2) at (1,-1) {};
	\node[shape=circle,fill=black,inner sep=1.5] (-1) at (2,0) {};
	\node[shape=circle,fill=black,inner sep=1.5] (0) at (-1,1) {};
	\node[shape=circle,fill=black,inner sep=1.5] (1) at (0,0) {};
	\node[shape=circle,fill=black,inner sep=1.5] (2) at (1,1) {};
	\node[shape=circle,fill=black,inner sep=1.5] (3) at (2,2) {};
	\node[shape=circle,fill=black,inner sep=1.5] (4) at (3,3) {};
	\node[shape=circle,fill=black,inner sep=1.5] (6) at (0,2) {};
	\node[shape=circle,fill=black,inner sep=1.5] (7) at (1,3) {};
	\node[shape=circle,fill=black,inner sep=1.5] (8) at (2,4) {};
	\node[shape=circle,fill=black,inner sep=1.5] (11) at (0,4) {};
	\node[shape=circle,fill=black,inner sep=1.5] (12) at (1,5) {};
	\node[shape=circle,fill=black,inner sep=1.5] (13) at (0,6) {};
	\node[shape=circle,fill=black,inner sep=1.5] (14) at (-1,7) {};
	\draw [thick] (-4) -- (-3) -- (-2) -- (-1) -- (2);
	\draw [thick] (2) -- (6) -- (0) -- (1) -- (-2);
	\draw [thick] (1) -- (2) -- (3) -- (4) -- (8) -- (12) -- (13) -- (14);
	\draw [thick] (3) -- (7) -- (6);
	\draw [thick] (7) -- (8);
	\draw [thick] (7) -- (11) -- (12);
	\node at (0,-5) {$\Lambda_{E_6}$};
\end{tikzpicture} \qquad  \vline \qquad \begin{tikzpicture}[scale=0.5]
	\node[shape=circle,fill=black,inner sep=1.5] (1) at (0,0) {};
	\node[shape=circle,fill=black,inner sep=1.5] (2) at (1,0) {};
	\node[shape=circle,fill=black,inner sep=1.5] (3) at (2,0) {};
	\node[shape=circle,fill=black,inner sep=1.5] (4) at (0.5,0.5) {};
	\node[shape=circle,fill=black,inner sep=1.5] (5) at (1.5,0.5) {};
	\node[shape=circle,fill=black,inner sep=1.5] (6) at (1,1) {};
	\node[shape=circle,fill=black,inner sep=1.5] (7) at (2,1) {};
	\node[shape=circle,fill=black,inner sep=1.5] (8) at (1.5,1.5) {};
	\node[shape=circle,fill=black,inner sep=1.5] (9) at (2.5,1.5) {};
	\node[shape=circle,fill=black,inner sep=1.5] (10) at (1,2) {};
	\node[shape=circle,fill=black,inner sep=1.5] (11) at (2,2) {};
	\node[shape=circle,fill=black,inner sep=1.5] (12) at (1.5,2.5) {};
	\node[shape=circle,fill=black,inner sep=1.5] (13) at (1,3) {};
	\node[shape=circle,fill=black,inner sep=1.5] (14) at (0.5,3.5) {};
	\node[shape=circle,fill=black,inner sep=1.5] (15) at (0,4) {};
	\node[shape=circle,fill=black,inner sep=1.5] (-4) at (0.5,-0.5) {};
	\node[shape=circle,fill=black,inner sep=1.5] (-5) at (1.5,-0.5) {};
	\node[shape=circle,fill=black,inner sep=1.5] (-6) at (1,-1) {};
	\node[shape=circle,fill=black,inner sep=1.5] (-7) at (2,-1) {};
	\node[shape=circle,fill=black,inner sep=1.5] (-8) at (1.5,-1.5) {};
	\node[shape=circle,fill=black,inner sep=1.5] (-9) at (2.5,-1.5) {};
	\node[shape=circle,fill=black,inner sep=1.5] (-10) at (1,-2) {};
	\node[shape=circle,fill=black,inner sep=1.5] (-11) at (2,-2) {};
	\node[shape=circle,fill=black,inner sep=1.5] (-12) at (1.5,-2.5) {};
	\node[shape=circle,fill=black,inner sep=1.5] (-13) at (1,-3) {};
	\node[shape=circle,fill=black,inner sep=1.5] (-14) at (0.5,-3.5) {};
	\node[shape=circle,fill=black,inner sep=1.5] (-15) at (0,-4) {};
	\draw [thick] (1) -- (4) -- (2) -- (5) -- (3);
	\draw [thick] (4) -- (6) -- (5) -- (7) -- (8) -- (6);
	\draw [thick] (7) -- (9) -- (11) -- (8) -- (10) -- (12) -- (11);
	\draw [thick] (12) -- (13) -- (14) -- (15);
	\draw [thick] (1) -- (-4) -- (2) -- (-5) -- (3);
	\draw [thick] (-4) -- (-6) -- (-5) -- (-7) -- (-8) -- (-6);
	\draw [thick] (-7) -- (-9) -- (-11) -- (-8) -- (-10) -- (-12) -- (-11);
	\draw [thick] (-12) -- (-13) -- (-14) -- (-15);
	\node at (1.5,-4.75) {$\Lambda_{E_7}$};
\end{tikzpicture}
\end{center}
\caption{The other minuscule posets beyond the rectangle and shifted staircase.} \label{fig:minuscule_posets}
\end{figure}
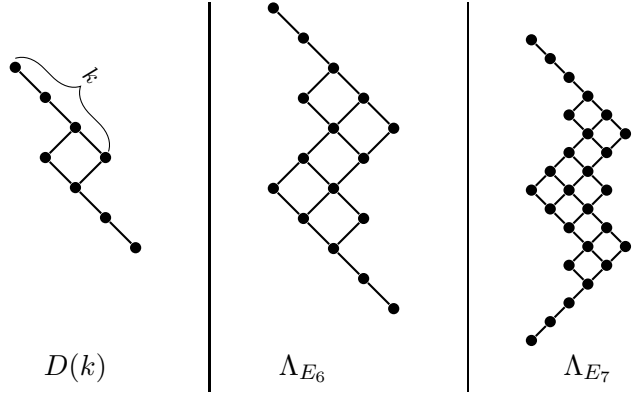

Evidently the rectangle and the shifted staircase are special. Their specialness extends to all the \dfn{minuscule posets}. The minuscule posets are a certain class of posets coming from the representation theory of simple Lie algebras, which exhibit remarkable combinatorial properties. We will not go over their formal definition. Instead, let us briefly review their classification. In addition to the rectangles and shifted staircases, there is one other infinite family of minuscule posets, the ``propeller poset'' $D(k)$ for $k\geq 2$, and there are also two exceptional minuscule posets $\Lambda_{E_6}$ and~$\Lambda_{E_7}$. These are depicted in \cref{fig:minuscule_posets}. Observe that all minuscule posets are graded and self-dual.

It was Proctor~\cite{proctor1984bruhat} who first recognized the combinatorial significance of the minuscule posets in the context of plane partition enumeration. We refer to that paper (see also~\cite{stembridge1994minuscule}) for further background on minuscule posets. Proctor proved that all the minuscule posets have product formulas for the $q$-enumeration of their bounded (reverse) $P$-partitions. That is, he showed that if $P$ is a minuscule poset then
\begin{equation} \label{eqn:rpp_minuscule}
\sum_{\pi \in \RPP_m(P)} q^{|\pi|} = \prod_{p \in P} \frac{[\rk(p)+m+1]_q}{[\rk(p)+1]_q}
\end{equation}
for any $m\geq 0$. Observe how \eqref{eqn:rpp_minuscule} includes \eqref{eqn:macmahon} and \eqref{eqn:rpp_shifted_staircase}. (In fact, Proctor conjectured that minuscule posets are the \emph{only} posets which have product formulas for $\sum_{\pi \in \RPP_m(P)} q^{|\pi|}$ of this form.) Notice that the limit $m\to \infty$ of \eqref{eqn:rpp_minuscule}, together with \eqref{eqn:rpp_gf}, implies that for a minuscule poset $P$ we also have
\begin{equation} \label{eqn:lin_minuscule}
\sum_{T\in\Lin(P)}q^{\comaj(T)} = [\#P]_q!\cdot \prod_{p \in P} \frac{1}{[\rk(p)+1]_q}.
\end{equation}

Once again, our barely set-valued results extend to all minuscule posets, because Defant et al.~proved:

\begin{thm}[{Defant et al.~\cite{defant2021homomesy}}] \label{thm:minuscule_ddeg_expression}
For $P$ a minuscule poset, there are constants $c(q) ,c_p(q) \in \Q(q)$ for which
\[ \ddeg = c(q) + \sum_{p\in P} c_p(q) \Tq{p}{q}.\]
\end{thm}

The case $q=1$ of \cref{thm:shifted_ddeg_expression} was again proved in \cite{hopkins2017CDE} (see also~\cite{rush2016minuscule}).

\begin{cor} \label{cor:minuscule_bsv_rpp}
For $P$ a minuscule poset, and any $m \geq 1$,
\[\sum_{\tau\in \RPPo_m(P)}q^{|\tau|-1} =   \frac{\sum_{p\in P} q^{\rk(p)} }{[\rk(P)+2]_q} \cdot [m]_q \cdot \sum_{\pi\in \RPP_m(P)}q^{|\pi|}.\] 
\end{cor}

\begin{cor} \label{cor:minuscule_bsv_lin}
For $P$ a minuscule poset,
\[\sum_{S\in \Lino(P)}q^{\comaj^{+1}(S)} = \frac{\sum_{p\in P} q^{\rk(p)} }{[\rk(P)+2]_q} \cdot [\#P+1]_q \cdot \sum_{T\in \Lin(P)}q^{\comaj(T)}.\] 
\end{cor}

\begin{proof}[Proofs of \cref{cor:minuscule_bsv_rpp,cor:minuscule_bsv_lin}]
As before, these follow from combining \cref{thm:minuscule_ddeg_expression} with \cref{cor:bsv_rpp_expectation,cor:bsv_lin_expectation}. To compute what the constant $c(q)$ must be, we use the strategy mentioned in~\cref{rem:rank_expectation} of computing $\E_{\murk{q}}(\ddeg)$. This computation gives $\E_{\murk{q}}(\ddeg) = \sum_{p\in P} q^{\rk(P)-\rk(p)}/[\rk(P)+2]_q$. But because every minuscule poset is self-dual, we have $\sum_{p\in P} q^{\rk(P)-\rk(p)} = \sum_{p\in P} q^{\rk(p)}$, which then yields the claimed formulas.
\end{proof}

As before, in light of~\eqref{eqn:rpp_minuscule} and~\eqref{eqn:lin_minuscule}, \cref{cor:minuscule_bsv_rpp,cor:minuscule_bsv_lin} give product formulas for $\sum_{\tau\in \RPPo_m(P)}q^{|\tau|-1}$ and $\sum_{S\in \Lino(P)}q^{\comaj^{+1}(S)}$ for any minuscule poset $P$.

\section{Final remarks} \label{sec:final}

We conclude with some final remarks and questions related to our work.

\begin{remark}
As mentioned in the introduction, the origin of all of this recent activity regarding product formulas for set-valued tableaux is in Brill--Noether theory~\cite{chan2018genera, chan2021euler}. We have no idea if our $q$-analogs have any significance in the context of Brill--Noether theory. It would be quite interesting if they did.
\end{remark}

\begin{remark}
Ultimately our results concern identities regarding certain polynomials in $q$. Hence it should be possible to remove all language concerning probability distributions, expectations of random variables, et cetera, from our arguments and work purely at the level of polynomials. Still, we find the probabilistic framing to be an important source of intuition and an aid in understanding our results.
\end{remark}

\begin{remark}
All of the results in this paper concern naturally labeled posets. However, the theory of $P$-partitions which underlies our work can be extended to the theory of $(P,\omega)$-partitions, where the labeling of the posets is no longer required to be natural (see \cite[Section 3.15]{stanley1996ec1}). It would be interesting to see which of our results can be extended to the more general setting of labeled posets. For example, it seems likely that a version of \cref{thm:m-weight} should hold for any labeling $\omega$. However, for general $\omega$ the corresponding probability distribution (at least, naively defined) will lack the crucial property of $q$-toggle-symmetry.
\end{remark}

\begin{remark}
Many of our proofs rely on results from~\cite{defant2021homomesy}, and especially on a  massive cancellation when computing the down-degree expectations that a priori has no bijective meaning. It is reasonable to ask if there exist bijective proofs of our formulas that do not rely on this cancellation. Relatedly, it would be interesting to find an efficient method for randomly sampling from $\mathcal{J}(P)$ according to the distributions $\murpp{q}{m}$ and~$\mulin{q}$.
\end{remark}

\begin{remark}
The results of~\cite{chan2017expected} (see also~\cite{reiner2018poset}) imply that
\[ \frac{\#\SYTo(\lambda)}{(n+1) \cdot \#\SYT(\lambda)} =\frac{ab}{a+b} \]
for all the so-called \dfn{``balanced''} partitions $\lambda$ of $n$, where $a$ is the number of rows of $\lambda$ and $b$ the number of columns of $\lambda$. (Note that, by \cref{cor:bsv_lin_expectation}, this ratio is $\E_{\mulin{1}}(\ddeg)$.) The balanced partitions are defined by a certain condition on their inner corners and include rectangles. Thus, this formula extends the product formula~\eqref{eqn:bsv_rect} for the number of barely set-valued SYTs of rectangular shape. Similarly, product formulas enumerating barely set-valued SYTs for so-called \dfn{``shifted balanced''} shapes were obtained in~\cite{hopkins2017CDE}, and for \dfn{shifted trapezoids} were obtained in~\cite{kim2020enumeration}. Hence one might wonder if there could be $q$-analogs for the enumeration of barely set-valued SYTs of other shapes beyond the ones we have considered here. We have not been able to find such $q$-analogs. 
\end{remark}

\begin{remark}
Some $q$-analogs of down-degree expectations of the kind we studied in this paper are proposed in~\cite[Appendix~A]{kim2020enumeration}. More precisely, the results of~\cite{chan2017expected} imply that
\[\frac{\#\RPPo_1(P_\lambda)}{\#\RPP_1(P_\lambda) }= \frac{ab}{a+b}\]
for balanced shapes $\lambda$ with $a$ rows and $b$ columns. (Note that, by \cref{cor:bsv_rpp_expectation}, this ratio is $\E_{\muuni{1}}(\ddeg)$.) In~\cite[Appendix~A]{kim2020enumeration} the authors propose a $q$-analog of this formula. It is not clear to us how their conjectural $q$-analogs are connected to our present work, although it would certainly be interesting to find a connection.
\end{remark}

\begin{remark}
In the direction of extending our work to other posets, we can ask for which posets $P$ there is a solution to
\[ \ddeg = c(q) + \sum_{p\in P}c_p(q) \Tq{p}{q}\]
with $c(q), c_p(q) \in \Q(q)$. We have checked by computer that there are no shapes beyond rectangles and no shifted shapes beyond shifted staircases for which a solution exists, among the shapes/shifted shapes with at most $20$ boxes. By contrast, when $q=1$ there are many more such shapes/shifted shapes (namely, the balanced/shifted balanced shapes of~\cite{chan2017expected,hopkins2017CDE}).
\end{remark}

\bibliography{q-barely-set-valued}{}
\bibliographystyle{abbrv}

\end{document}